\documentclass[review]{elsarticle}

\usepackage{amsmath}
\usepackage{amssymb}
\usepackage{mathrsfs}
\usepackage{graphicx}
\usepackage{amsmath}
\usepackage{latexsym}
\usepackage{amsfonts}
\usepackage{multirow}
\usepackage{cases}
\usepackage{tabularx}
\usepackage{booktabs}
\usepackage{color}
\usepackage{setspace}
\usepackage{algorithm}
\usepackage{algorithmic}

\newtheorem{theorem}{Theorem}
\newtheorem{definition}{Definition}
\newtheorem{proposition}{Proposition}
\newtheorem{lemma}{Lemma}

\newtheorem{corollary}{Corollary}

\newcommand{\kgg}[1]{\mbox{}\hspace{#1}}
\newcommand{\jbe}{\begin{equation}}
\newcommand{\jen}{\end{equation}}

\newcommand{\kg}{\mbox{}\hspace{0.3in}}

\begin{document}
\begin{frontmatter}

\title{Optimization Condition and Algorithm of Optimization with Convertible Nonconvex Function\tnoteref{mytitlenote}}
\tnotetext[mytitlenote]{This work was supported by the National Natural Science
Foundation of China (Grant No. 11871434).}

%


\author[1]{Min Jiang}
\author[2]{Rui Shen}
\author[1]{Zhiqing Meng}
\author[3]{Chuangyin Dang}
\address[1]{School of Management, Zhejiang University of Technology, Hangzhou, Zhejiang, 310023, China}
\address[2]{School of Economics, Zhejiang University of Technology, Hangzhou, Zhejiang, 310023, China}
\address[3]{Department of Advanced Design and System Engineering, City University of Hong Kong, Kowloon, Hong Kong}
\date{Jan. 9, 2022}

\begin{abstract}
The paper introduces several new concepts for solving  nonconvex or nonsmooth optimization problems, including convertible nonconvex function, exact convertible nonconvex function and differentiable convertible nonconvex function. It is proved herein many nonconvex functions or nonsmooth (or discontinuous) functions are actually  convertible nonconvex functions and  convertible nonconvex function operations such as addition, subtraction, multiplication or division result in convertible nonconvex functions.  The sufficient condition for judging a global optimal solution to unconstrained optimization problems with differentiable convertible nonconvex functions is proved, which is equivalent to  Karush-Kuhn-Tucker(KKT) condition. Two Lagrange functions of differentiable convertible nonconvex function are defined with their dual problems  defined accordingly. The strong duality theorem is proved, showing that the optimal objective value of the global optimal solution is equal to the optimal objective value of the dual problem, which is equivalent to KKT condition. An augmented Lagrangian penalty function algorithm is proposed and its convergence is proved. So the paper provides a new idea for solving unconstrained nonconvex or non-smooth optimization problems and avoids  subdifferentiation or smoothing techniques by using some gradient search algorithms, such as gradient descent algorithm, Newton algorithm and so on.
\end{abstract}

\begin{keyword} convertible nonconvex function\sep convertible nonconvex function optimization\sep global optimal solution\sep sufficient condition\sep Lagrange dual\sep algorithm.
\end{keyword}

\end{frontmatter}

\begin{abstract}

\end{abstract}

\section{Introduction}
The following optimization problem  is defined by
\begin{eqnarray*}
\mbox{(P1)} &\min & f(\boldsymbol{x})\\
&s.t.& \boldsymbol{x}\in X\subset R^n.
\end{eqnarray*}
When $X=R^n$, (P1) is unconstrained optimization.  When $f$ is a nonconvex function, it is difficult to find  a global optimal solution to (P1).  In 1976,  McCormick  presented a method  to find a global optimal solution to (P1) when $f$ is a factorable function, where
(P1) is called factorable programming\cite{McCormick}. The main idea of factorable programming is that when a convergent sequence of global solutions $\boldsymbol{x}^k$ to a sequence of convex underestimating programming is obtained, $\boldsymbol{x}^*$ is a global optimal solution to (P1) as $\boldsymbol{x}^k\to \boldsymbol{x}^*$, as shown in \cite{Bunin,He,McCormick,Serrano}. Hence, Sherali and Wang(2001)\cite{Sherali}, Tawarmalani and Sahinidis(2004)\cite{Tawarmalani}, and Granvilliers(2020)\cite{Granvilliers} showed a branch and bound algorithm for solving a global optimal solution to factorable programming. In 1988, Jackson and Mccormick proved the results of second-order sensitivity analysis of factorable programming guarantee the stability of the optimal solution  to factorable programming\cite{Jackson1}. Then, factorable programming was applied in encyclopedia of operations research and management science \cite{Jackson2}. Furthermore, Hascot et al.(2013)\cite{Hascot} proved that factorable programming is suitable for designing an algorithm to solve large-scale optimization problems. All the above shows that there are effective algorithms for solving a global solution to factorable programming.

Hanson and Mond(1987)\cite{Hanson} pointed out that only the nonconvex optimization with special structure can find a global optimal solution to (P1), i.e. when $f$ is not convex but convex transformable (so (P1) is called convex transformable programming). Since convex transformable programming is converted into a relaxed convex optimization, Hirschberger(2005)\cite{Hirschberger},  Neto et al.(2006)\cite{Neto},  Andreas et al.(2009)\cite{Andreas}, Khajavirad et al.(2014)\cite{Khajavirad} and  Nohra and Sahinidis(2018)\cite{Nohra} gave  algorithms to solve a global solution. 
Grant et al.(2006)\cite{Grant} gave out  methods to solve nonconvex optimization problems, and Shen et al.(2017)\cite{Shen} gave a multiconvex programming which is a nonconvex optimization problem composed of multiconvex functions.
In recent years, in  engineering, various practical models of  multiconvex programming are studied as shown in \cite{Chiu,Hours,Ichihara,Suh}. But, to find a global optimal solution to multiconvex programming is  more difficult than to find one to other nonconvex optimization problems.

In summary, some convex transformable functions are factorable as shown in \cite{Khajavirad}.  The above special nonconvex functions have some weaknesses. For example, they require a complex technique to construct convex underestimating programming of factor programming, which makes it difficult to construct its global algorithm. Besides, Khajavirad et al pointed our that it is difficult to judge whether a nonconvex function is convex transformable \cite{Khajavirad}. In addition, the major obstacle is that the optimal condition of the global optimal solution to nonconvex and nondifferentiable optimization problem (P1) has not been found in theory.

 On the other hand, subdifferential theory and smoothing techniques are widely used to solve nonconvex and nondifferentiable optimization problem in \cite{Chen1,Chen2,Clarke,Huang}. When $f$ in (P1) is nonconvex or nonsmooth, it is difficult to find its a global solution under optimization condition or dual problem.

In this paper, if $f(\boldsymbol{x})$ is nonconvex, nonsmooth or discontinuous on $\boldsymbol{x}\in R^n$, an idea of mapping  $f$ from $R^n$ to a higher dimensional space $R^n\times R^m$ such that $f(\boldsymbol{x})=g(\boldsymbol{x},\boldsymbol{y})$ with $(\boldsymbol{x},\boldsymbol{y})\in X(g)\subset R^n\times R^m$ in Definition 2.1, where
 $g$ is a differentiable and convex, is proposed to solve (P1). So,   $$\min\limits_{\boldsymbol{x}\in X} f(\boldsymbol{x})=\min\limits_{(\boldsymbol{x},\boldsymbol{y})\in X(g)} g(\boldsymbol{x},\boldsymbol{y}). $$ When $X(g) $ is expected to have a good structure,  a global or local optimal solution to  $\min\limits_{(\boldsymbol{x},\boldsymbol{y})\in X(g)} g(\boldsymbol{x},\boldsymbol{y})$ may be found. The paper applies this idea and proves theoretically that the sufficient conditions and Lagrangian dual theorem  can be realized for solving a global or local optimal solution to the unconstrained optimization (P1) to avoid using subdifferentiation or smoothing techniques, which allows us to directly use some gradient search algorithms, such as gradient descent algorithm, Newton algorithm and so on.

So, the main contributions include the three results as follows. (1) The definition of convertible nonconvex function, exact convertible nonconvex function and differentiable convertible nonconvex function, which covers a large scale of nonconvex or discontinuous functions, such as convex transformable functions, multiconvex functions, 0-norm functions, etc. (2) The sufficient conditions and the duality theorem of a global optimal solution to  (P1) are proved in theory. However, it is impossible to obtain these results if the convertible function form technology is not used. (3) An augmented Lagrangian penalty function algorithm is designed to find an approximately optimal solution to (P1) where some gradient search algorithms are directly used.

The remainder of this paper is organized as follows. In Section 2, the definition, some examples and some properties of (exact, differentiable) convertible nonconvex function are given. In Section 3, unconstrained optimization with a differentiable convertible nonconvex function is defined and the sufficient conditions of its global optimal solution are proved. In Section 4, the  Lagrange function of a differentiable convertible nonconvex function and  dual problems of unconstrained optimization are defined. Their dual properties are discussed. An augmented Lagrangian penalty function algorithm and its convergence are given. In Section 5, a conclusion is given.

\section{Convertible Nonconvex Function}

Convertible nonconvex function, exact convertible nonconvex function and differentiable convertible nonconvex function are defined first. Then examples are given to show some nonconvex, nonsmooth or discontinuous functions are differentiable  exact convertible nonconvex ones.

\begin{definition}
 Let $s+r+1$ functions: $g, g_i,h_j:R^n\times R^m\to R^1$($i=1,2,\cdots,s;j=1,2,\cdots,r$) be convex and $f: R^n \rightarrow R$ be nonconvex or nonsmooth. Let a set
\begin{eqnarray}
X(f)=\{(\boldsymbol{x},\boldsymbol{y})| g_i(\boldsymbol{x},\boldsymbol{y})\leq 0,i=1,2,\cdots,s; h_j(\boldsymbol{x},\boldsymbol{y})=0, j=1,2,\cdots,r\}\label{eq:d1}
\end{eqnarray}
For a fixed $\boldsymbol{x}$, let
\begin{eqnarray}
X(f(\boldsymbol{x}))=\{\boldsymbol{y}\in R^m|(\boldsymbol{x},\boldsymbol{y})\in X(f)\}\label{eq:d2}
\end{eqnarray}
and
\begin{eqnarray}
X(g)=\{(\boldsymbol{x},\boldsymbol{y})\in X(f)|g(\boldsymbol{x},\boldsymbol{y})=f(\boldsymbol{x})\}.\label{eq:d3}
\end{eqnarray}

$\bullet$  If for each $(\boldsymbol{x},f(\boldsymbol{x}))$ there is a $\bar{\boldsymbol{y}}\in X(f(\boldsymbol{x}))$ such that $f(\boldsymbol{x})=g(\boldsymbol{x},\bar{\boldsymbol{y}})$ and
\begin{eqnarray*}
f(\boldsymbol{x})=\min\limits_{\boldsymbol{y}\in X(f(\boldsymbol{x}))} g(\boldsymbol{x},\boldsymbol{y}),
\end{eqnarray*}
then  $f$ is called a convertible nonconvex function. $[g:g_1,g_2,\cdots,g_s;h_1,h_2,\cdots,h_r]$ is  convertible nonconvex form of $f$, briefing as $f=[g:g_1,g_2,\cdots,g_s;h_1,h_2,\cdots,h_r]$. Moreover, the following concepts are defined.

$\bullet$ If $f$ is continuous on $R^n$ and $g,g_1,g_2,\cdots,g_s,h_1,h_2,\cdots,h_r$ are continuous on $R^n\times R^m$,  $f$ is called strong continuous convertible nonconvex function.

$\bullet$ If $s+r+1$ functions: $g, g_i,h_j$($i=1,2,\cdots,s;j=1,2,\cdots,r;s>0$) are differentiable,
  $f$ is called a differentiable convertible nonconvex function.

$\bullet$ If $f$ is a differentiable convertible nonconvex function and $f(\cdot)$ is a differentiable function on $R^n$,  $f$ is  a strong differentiable convertible nonconvex function.

$\bullet$ If $g(\boldsymbol{x},\boldsymbol{y})=f(\boldsymbol{x})$ for any $(\boldsymbol{x},\boldsymbol{y})\in X(f)$,  then $f$ is called an exact convertible nonconvex function and $(\boldsymbol{x},\boldsymbol{y})$ is called a mappable point of $f$.
\end{definition}

{\bf Remark } The definition of  convertible nonconvex function differs form that of upper(lower)-$C^k$ function  \cite{Daniilidis}, where the upper(lower)-$C^k$ function is continuous. For example, 0-norm $\|\boldsymbol{x}\|_0$ is not upper-$C^k$ function \cite{Daniilidis}, because it is not continuous. But, $\|\boldsymbol{x}\|_0$ is a convertible nonconvex function in Example 4.

In this paper, we always assume that $f$ is a convertible nonconvex function.
If $s=0$, inequality does not exist in the convertible nonconvex form of $f$.
If $f$ is an exact convertible nonconvex function, $X(f)=X(g)$.

Some examples of convertible nonconvex functions are given as follows.

{\bf Example 1} Two convertible nonconvex forms of $f(x_1,x_2)=|x_1x_2|^\frac{1}{3}+x_1^2+x_2^2$ are defined  by
\begin{eqnarray*}
f=[g:g_1;h_1,h_2,h_3,h_4]=[y_4+x_1^2+x_2^2:
-y_4;0.5(x_1+x_2)^2-y_1-0.5y_2,\\ x_1^2+x_2^2-y_2,y_1^2-y_3,y_4^6-y_3]
\end{eqnarray*}
and
\begin{eqnarray*}
f=[g:g_1;h_1,h_2,h_3,h_4]=[y_4+x_1^2+x_2^2:
-y_4;0.25(x_1+x_2)^2-y_1-0.25y_2,\\ (x_1-x_2)^2-y_2,y_1^2-y_3,y_4^6-y_3]
\end{eqnarray*}
respectively.  Then, $f=[g:g_1;h_1,h_2,h_3,h_4]$ is differentiable exact convertible nonconvex. And, the convertible nonconvex form of  $f$ is not only one. $f$ is noncovex and nonsmooth.


{\bf Example 2} Let $\boldsymbol{b}_1,\boldsymbol{b}_2\in R^n$ be given. A convertible nonconvex forms of $f(\boldsymbol{x})=(\sqrt{|\boldsymbol{b}_1^\top \boldsymbol{x})|}-\sqrt{|\boldsymbol{b}_2^\top \boldsymbol{x}|})^2$ is defined by
\begin{eqnarray*}
f=[g:g_1,g_2;h_1,h_2,h_3,h_4]=(y_1-y_2)^2: -y_1,
-y_2;y_1^2-y_3,y_2^2-y_4,\\
(\boldsymbol{b}_1^\top \boldsymbol{x})^2-y_3,(\boldsymbol{b}_2^\top \boldsymbol{x})^2-y_4].
\end{eqnarray*}
$f$ is differentiable exact convertible nonconvex.

{\bf Example 3}  Let multi-classification function $f(\boldsymbol{x})=\sum\limits_{i=1}^I(|\boldsymbol{a}_i^\top\boldsymbol{x}|-b_i)^2$ be nonsmooth and nonconvex in \cite{Mohri}, where sample data $(\boldsymbol{a}_i,b_i)\in R^n\times R^1(i=1,2,\cdots,I)$ are given. Let
\begin{eqnarray*}
f(\boldsymbol{x})=g(\boldsymbol{x},\boldsymbol{y})&=&\sum\limits_{i=1}^I(y_i-b_i)^2: \\
g_{i}(\boldsymbol{x},\boldsymbol{y})&=&-y_i\leq 0,\ i=1,2,\cdots,I;\\
h_{i}(\boldsymbol{x},\boldsymbol{y})&=&y_i^2-y_{i+I}=0,\ i=1,2,\cdots,I,\\
h_{i+I}(\boldsymbol{x},\boldsymbol{y})&=&(\boldsymbol{a}_i^\top\boldsymbol{x})^2-y_{i+I}=0, \ i=1,2,\cdots,I,
\end{eqnarray*}
where $\boldsymbol{y}\in R^{2I}$. $f=[g:g_1,\cdots,g_I;h_1,\cdots,h_{2I}]$ is  twice differentiable exact  convertible nonconvex.

We have the following properties.

\begin{proposition}
If $f,f_1,f_2:R^n\to R^1$ are convertible nonconvex, then the following conclusions hold,

(i) $\alpha f_1$ is convertible nonconvex, where $\alpha\not=0$.

(ii) $\alpha_1 f_1+\alpha_2 f_2$ is convertible nonconvex, where $\alpha_1,\alpha_2 \not=0$.

(iii) $f_1f_2$ is convertible nonconvex.

(iv) $\phi(f)$ is convertible nonconvex, where $\phi:R^1\to R^1$ is a monotone increasing convex function.

(v) $\max\{f_1,f_2\}$ is convertible nonconvex.

(vi)  $\frac{f_1}{f_2}$ is convertible nonconvex.

(vii) $-log(f_1)$ is convertible nonconvex.

(viii) $exp(af_1)$ is convertible nonconvex, where $a>1$.
\end{proposition}

By Proposition 1, all polynomial functions are convertible nonconvex. By Proposition 1, G-convex functions in Example 1-7 in \cite{Khajavirad} are convertible nonconvex. Some multi-convex functions are  convertible nonconvex, such as $f(\boldsymbol{x})=x_1x_2\cdots x_n$.  Therefore, convertible nonconvex functions include a wide range of non-convex functions.  By Proposition 2.1, the following corollaries hold.

\begin{corollary}
If $f,f_1,f_2:R^n\to R$ are exact convertible nonconvex, then the results of all operations in  Proposition 2.1 are exact convertible nonconvex.
\end{corollary}



In the following example, the nonsmooth, nonconvex and discontinuous function is a differentiable convertible nonconvex function.

{\bf Example 4}  Let 0-norm function $f(\boldsymbol{x})=\lambda\|\boldsymbol{x}\|_0+b(\boldsymbol{x})$ in machine learning be nonsmooth nonconvex and discontinuous in \cite{Chen1,Chen2}, where $b:R^n\to R^1$ is convex. There are $\boldsymbol{y} \in R^{2n}$
such that
\begin{eqnarray*}
f(\boldsymbol{x})=g(\boldsymbol{x},\boldsymbol{y})&=&\lambda\sum\limits_{i=1}^n y_i+b(\boldsymbol{x}):\\
g_{i}(\boldsymbol{x},\boldsymbol{y})&=&-y_i\leq 0,i=1,2,\cdots,n, \\
g_{i+n}(\boldsymbol{x},\boldsymbol{y})&=&y_i\leq 1, i=1,2,\cdots,n;\\
h_j(\boldsymbol{x},\boldsymbol{y})&=&(x_i+y_i-1)^2-y_{i+n}=0,i=1,2,\cdots,n,\\
h_{i+n}(\boldsymbol{x},\boldsymbol{y})&=&x_i^2+(y_i-1)^2-y_{i+n}=0,i=1,2,\cdots,n,\\
h_{i+2n}(\boldsymbol{x},\boldsymbol{y})&=&y_i^2-y_i=0,i=1,2,\cdots,n.
\end{eqnarray*}
So, $f$ is twice differentiable convertible nonconvex. But, $f$ is not exact convertible nonconvex.

The above examples tell us that some nonsmooth, nonconvex or discontinuous functions may be twice differentiable convertible nonconvex functions such that many of nonsmooth ,nonconvex optimization problems can be converted to differentiable convertible nonconvex optimization problems. Hence, it is very  meaningful to study differentiable convertible nonconvex optimization problems.

\section{Unconstrained Optimization with Convertible Nonconvex Function}

In this section, it is always assumed that $f$ is a differentiable convertible nonconvex function with $f=[g:g_1,g_2,\cdots; g_s, h_1,h_2,\cdots,h_r]$.
But, $f(\boldsymbol{x})$ is not necessarily differentiable or convex.

The following unconstrained optimization problem  is considered:
\begin{eqnarray*}
\mbox{(UOP)}\qquad & \min\; & f(\boldsymbol{x}) \\
& \mbox{s.t.}\; & \boldsymbol{x}\in R^n,
\end{eqnarray*}
where $f: R^n \rightarrow R$ is  differentiable convertible nonconvex.  The constrained optimization problem of (CNP) is defined by
\begin{eqnarray*}
\mbox{(CNP)}\qquad & \min\; & g(\boldsymbol{x},\boldsymbol{y}) \\
& \mbox{s.t.}\; &  g_i(\boldsymbol{x},\boldsymbol{y})\leq 0,i=1,2,\cdots,s,\\
&& h_j(\boldsymbol{x},\boldsymbol{y})=0,j=1,2,\cdots,r,\\
&& \boldsymbol{x}\in R^n,\boldsymbol{y}\in R^m.
\end{eqnarray*}
(UOP) or (CNP) is called convertible nonconvex programming.
Since $f: R^n \rightarrow R$ is  convertible nonconvex, it is clear
$$\min\limits_{\boldsymbol{x}\in R^n} \ f(\boldsymbol{x})=\min\limits_{(\boldsymbol{x},\boldsymbol{y})\in X(f)} \ g(\boldsymbol{x},\boldsymbol{y})=\min\limits_{(\boldsymbol{x},\boldsymbol{y})\in X(g)} \ g(\boldsymbol{x},\boldsymbol{y}).$$
When $r=0$ or all $h_j(\boldsymbol{x},\boldsymbol{y})=0(j=1,2,\cdots,r)$ are linear, (CNP) is a convex programming.

Let $\boldsymbol{d}=(\boldsymbol{d}_1,\boldsymbol{d}_2)\in R^n\times R^m$. The linear programming of (CNP) at a fixed $(\boldsymbol{x},\boldsymbol{y})$ is defined by
\begin{eqnarray*}
\mbox{(CNP)}(\boldsymbol{x},\boldsymbol{y})\qquad & \min\; & \nabla g(\boldsymbol{x},\boldsymbol{y})^\top \boldsymbol{d} \\
& \mbox{s.t.}\; & g_i(\boldsymbol{x},\boldsymbol{y})+\nabla  g_i(\boldsymbol{x},\boldsymbol{y})^\top\boldsymbol{d}\leq 0,i=1,2,\cdots,s,\\
&& \nabla h_j(\boldsymbol{x},\boldsymbol{y})^\top \boldsymbol{d}\leq 0,j=1,2,\cdots,r,\\
         &&      \boldsymbol{d}\in R^n\times R^m.
\end{eqnarray*}

Let set $I=\{1,2,\cdots,s\}$, $J=\{1,2,\cdots,r\}$ and vector functions
\begin{eqnarray*}
\boldsymbol{g}(\boldsymbol{x},\boldsymbol{y})=(g_1(\boldsymbol{x},\boldsymbol{y}),g_2(\boldsymbol{x},\boldsymbol{y}),\cdots,g_s(\boldsymbol{x},\boldsymbol{y}))^\top
\end{eqnarray*}
and
\begin{eqnarray*}
\boldsymbol{h}(\boldsymbol{x},\boldsymbol{y})=(h_1(\boldsymbol{x},\boldsymbol{y}),h_2(\boldsymbol{x},\boldsymbol{y}),\cdots,h_r(\boldsymbol{x},\boldsymbol{y}))^\top.
\end{eqnarray*}
When $f(\boldsymbol{x})$ is  differentiable on $\boldsymbol{x}\in R^n$, it is easy to determine that  $\boldsymbol{x}$ is not the optimal solution to (UOP) when $\nabla f(\boldsymbol{x})\not=0$. However, it is very difficult to judge the global optimal solution to (UOP).
How to judge the global optimal solution of (CNP) is proved in the following.

\begin{theorem}
Suppose that $(\boldsymbol{x}^*,\boldsymbol{y}^*)\in X(g)$.
If there is an optimal solution $\boldsymbol{d}^*$ to (CNP)$(\boldsymbol{x}^*,\boldsymbol{y}^*)$ such that $\nabla g(\boldsymbol{x}^*,\boldsymbol{y}^*)^\top \boldsymbol{d}^* \geq 0$, then  $(\boldsymbol{x}^*,\boldsymbol{y}^*)$ is an optimal solution to (CNP), $\boldsymbol{x}^*$ is an optimal solution to (UOP) and there are $\boldsymbol{u}^*=(u_1^*,u_2^*,\cdots,u_s^*)^\top\geq 0$ and $\boldsymbol{v}^*=(v_1^*,v_2^*,\cdots,v_r^*)^\top\geq 0$ such that
\begin{eqnarray}
\nabla g(\boldsymbol{x}^*,\boldsymbol{y}^*)+ \nabla \boldsymbol{g}(\boldsymbol{x}^*,\boldsymbol{y}^*)^\top\boldsymbol{u}^* +\nabla \boldsymbol{h}(\boldsymbol{x}^*,\boldsymbol{y}^*)^\top\boldsymbol{v}^*=0.\label{eq:d4}\\
u^*_ig_i(\boldsymbol{x}^*,\boldsymbol{y}^*)^\top=0,i=1,2,\cdots,s.\label{eq:d5}
\end{eqnarray}
\end{theorem}

{\it Proof.} For any $(\boldsymbol{x},\boldsymbol{y})\in X(g)$, we have
\begin{eqnarray*}
g(\boldsymbol{x},\boldsymbol{y})-g(\boldsymbol{x}^*,\boldsymbol{y}^*)\geq \nabla g(\boldsymbol{x}^*,\boldsymbol{y}^*)^\top[(\boldsymbol{x},\boldsymbol{y})-(\boldsymbol{x}^*,\boldsymbol{y}^*)],\\
 0\geq g_i(\boldsymbol{x},\boldsymbol{y})\geq g_i(\boldsymbol{x}^*,\boldsymbol{y}^*)+\nabla g_i(\boldsymbol{x}^*,\boldsymbol{y}^*)^\top[(\boldsymbol{x},\boldsymbol{y})-(\boldsymbol{x}^*,\boldsymbol{y}^*)],\  i=1,2,\cdots,s,\\
 0=h_j(\boldsymbol{x},\boldsymbol{y})-h_j(\boldsymbol{x}^*,\boldsymbol{y}^*)\geq \nabla h_j(\boldsymbol{x}^*,\boldsymbol{y}^*)^\top[(\boldsymbol{x},\boldsymbol{y})-(\boldsymbol{x}^*,\boldsymbol{y}^*)],\ j=1,2,\cdots,r.
\end{eqnarray*}
So, $(\boldsymbol{x},\boldsymbol{y})-(\boldsymbol{x}^*,\boldsymbol{y}^*)$ is a feasible solution to (CNP)$(\boldsymbol{x}^*,\boldsymbol{y}^*)$, then
\begin{eqnarray*}
f(\boldsymbol{x})-f(\boldsymbol{x}^*)&=&g(\boldsymbol{x},\boldsymbol{y})-g(\boldsymbol{x}^*,\boldsymbol{y}^*)\\
&\geq& \nabla g(\boldsymbol{x}^*,\boldsymbol{y}^*)^\top[(\boldsymbol{x},\boldsymbol{y})-(\boldsymbol{x}^*,\boldsymbol{y}^*)]\\
&\geq& \nabla g(\boldsymbol{x}^*,\boldsymbol{y}^*)^\top\boldsymbol{d}^*\geq 0.
\end{eqnarray*}
Hence, $(\boldsymbol{x}^*,\boldsymbol{y}^*)$ is an optimal solution to (CNP) and $\boldsymbol{x}^*$ is an optimal solution to (UOP). Because (CNP)$(\boldsymbol{x}^*,\boldsymbol{y}^*)$ is linear programming, it is rewritten as
\begin{eqnarray*}
\mbox{(CNP)}(\boldsymbol{x}^*,\boldsymbol{y}^*)\qquad & \max\; & -\nabla g(\boldsymbol{x}^*,\boldsymbol{y}^*)^\top \boldsymbol{d} \\
& \mbox{s.t.}\; &  \nabla g_i(\boldsymbol{x}^*,\boldsymbol{y}^*)^\top\boldsymbol{d}\leq -g_i(\boldsymbol{x}^*,\boldsymbol{y}^*),\  i=1,2,\cdots,s,\\
&&\nabla h_j(\boldsymbol{x}^*,\boldsymbol{y}^*)^\top \boldsymbol{d}\leq 0,j=1,2,\cdots,r,\\
         &&      \boldsymbol{d}\in R^n\times R^m.
\end{eqnarray*}
So, the dual problem (DTCP)$(\boldsymbol{x}^*,\boldsymbol{y}^*)$ of (CNP)$(\boldsymbol{x}^*,\boldsymbol{y}^*)$  is defined as follows.
\begin{eqnarray*}
\mbox{(DTCP)}(\boldsymbol{x}^*,\boldsymbol{y}^*)\qquad & \min\; & - \sum\limits_{i=1}^s  u_i g_i(\boldsymbol{x}^*,\boldsymbol{y}^*)+ \sum\limits_{j=1}^s 0\cdot v_j \\
& \mbox{s.t.}\; &  \sum\limits_{i=1}^s u_i \nabla g_i(\boldsymbol{x}^*,\boldsymbol{y}^*)+\sum\limits_{j=1}^r v_i \nabla h_j(\boldsymbol{x}^*,\boldsymbol{y}^*) = -\nabla g(\boldsymbol{x}^*,\boldsymbol{y}^*)\\
         &&      u_i, v_j\geq 0, i=1,2,\cdots,s,j=1,2,\cdots,r;
\end{eqnarray*}
where $(u_1,u_2,\cdots,u_s)$ and $(v_1,v_2,\cdots,v_r)$  are dual variables. By the strong dual theorem of linear programming, there is an optimal solution $\boldsymbol{u}^*=(u_1^*,u_2^*,\cdots,u_r^*)^\top\geq 0$ and $\boldsymbol{v}^*=(v_1^*,v_2^*,\cdots,v_r^*)^\top\geq 0$ to $\mbox{(DTCP)}(\boldsymbol{x}^*,\boldsymbol{y}^*)$ such that
$$0 \leq - \sum\limits_{i=1}^s  u_i^* g_i(\boldsymbol{x}^*,\boldsymbol{y}^*)=-\nabla g(\boldsymbol{x}^*,\boldsymbol{y}^*)^\top \boldsymbol{d}^* \leq 0.$$
Hence, \eqref{eq:d4} and \eqref{eq:d5} are true. \\

So, \eqref{eq:d4} and \eqref{eq:d5} are a KKT condition of (CNP) of (UOP). Theorem 1 means that there is optimization condition of (UOP) if $f$ is not differentiable. Usually, if $f$ has subdifferentiation at optimal point $\boldsymbol{x}^*$, there may be an optimality condition $0\in \partial f(\boldsymbol{x}^*)$ in \cite{Clarke}.
 By Theorem 1, the following corollary is true.

\begin{corollary}
Suppose that $(\boldsymbol{x}^*,\boldsymbol{y}^*)\in X(g)$. Let the problem $\mbox{(CNPP)}(\boldsymbol{x}^*,\boldsymbol{y}^*)$
\begin{eqnarray*}
 & \min\; & \nabla g(\boldsymbol{x}^*,\boldsymbol{y}^*)^\top [(\boldsymbol{x},\boldsymbol{y})-(\boldsymbol{x}^*,\boldsymbol{y}^*)]\\
& \mbox{s.t.}\; &  g_i(\boldsymbol{x}^*,\boldsymbol{y}^*)+\nabla g_i(\boldsymbol{x}^*,\boldsymbol{y}^*)^\top[(\boldsymbol{x},\boldsymbol{y})-(\boldsymbol{x}^*,\boldsymbol{y}^*)]\leq 0,\  i=1,2,\cdots,s,\\
&&\nabla h_j(\boldsymbol{x}^*,\boldsymbol{y}^*)^\top  [(\boldsymbol{x},\boldsymbol{y})-(\boldsymbol{x}^*,\boldsymbol{y}^*)]\leq 0,j=1,2,\cdots,r,\\
         &&      (\boldsymbol{x},\boldsymbol{y})\in X(g).
\end{eqnarray*}
  If $(\boldsymbol{x}^*,\boldsymbol{y}^*)$  is an optimal solution  to (CNPP)$(\boldsymbol{x}^*,\boldsymbol{y}^*)$, then $(\boldsymbol{x}^*,\boldsymbol{y}^*)$  is an optimal solution  to (CNP) and $\boldsymbol{x}^*$ is an optimal solution to (UOP), i.e. if $(\boldsymbol{x}^*,\boldsymbol{y}^*)$  is not an optimal solution  to (CNP), then there is an $(\boldsymbol{x},\boldsymbol{y})\in X(g)$ such that
  $\nabla g(\boldsymbol{x}^*,\boldsymbol{y}^*)^\top [(\boldsymbol{x},\boldsymbol{y})-(\boldsymbol{x}^*,\boldsymbol{y}^*)]<0$.
\end{corollary}

 Because (CNPP)$(\boldsymbol{x}^*,\boldsymbol{y}^*)$ is not linear programming, it is  difficult  to solve it. The inverse proposition of Theorem 1 holds   as follows.

\begin{theorem}
 Let $(\boldsymbol{x}^*,\boldsymbol{y}^*)\in X(g)$.
If there are $\boldsymbol{u}^*=(u_1^*,u_2^*,\cdots,u_r^*)^\top\geq 0$ and $\boldsymbol{v}^*=(v_1^*,v_2^*,\cdots,v_r^*)^\top\geq 0$ such that \eqref{eq:d4}  and \eqref{eq:d5} hold,
then $(\boldsymbol{x}^*,\boldsymbol{y}^*)$ is an optimal solution to (CNP) and $\boldsymbol{x}^*$ is an optimal solution to (UOP).
\end{theorem}

The following conclusion is clear.

\begin{theorem}
Let $(\boldsymbol{x}^*,\boldsymbol{y}^*)\in X(g)$.
If $\nabla g(\boldsymbol{x}^*,\boldsymbol{y}^*)=0$, then  $(\boldsymbol{x}^*,\boldsymbol{y}^*)$ is an optimal solution to (CNP) and  $\boldsymbol{x}^*$ is an optimal solution to (UOP).
\end{theorem}


If the condition of Theorem 1 or Theorem 2 does not hold, the global optimal solution to (CNP) is judged by solving linear programming (CNP0)$(\boldsymbol{x},\boldsymbol{y})$.

The linear programming (CNP0)$(\boldsymbol{x},\boldsymbol{y})$ of (CNP) at a fixed $(\boldsymbol{x},\boldsymbol{y})$ is defined as follows.
\begin{eqnarray*}
\mbox{(CNP0)}(\boldsymbol{x},\boldsymbol{y})\qquad & \min\; & \nabla g(\boldsymbol{x},\boldsymbol{y})^\top \boldsymbol{d} \\
& \mbox{s.t.}\; & g_i(\boldsymbol{x},\boldsymbol{y})+\nabla  g_i(\boldsymbol{x},\boldsymbol{y})^\top\boldsymbol{d}\leq 0,i=1,2,\cdots,s,\\
&& \nabla h_j(\boldsymbol{x},\boldsymbol{y})^\top \boldsymbol{d}= 0,j=1,2,\cdots,r,\\
         &&      \boldsymbol{d}\in R^n\times R^m.
\end{eqnarray*}

Let two directional sets at a fixed $(\boldsymbol{x},\boldsymbol{y})$ be defined respectively by
\begin{eqnarray}
T(\boldsymbol{x},\boldsymbol{y})=\{\boldsymbol{d}\in R^n\times R^m&|& g_i(\boldsymbol{x},\boldsymbol{y})+\nabla  g_i(\boldsymbol{x},\boldsymbol{y})^\top\boldsymbol{d}\leq 0,i=1,2,\cdots,s,\nonumber\\
&&\nabla h_j(\boldsymbol{x},\boldsymbol{y})^\top \boldsymbol{d}\leq 0,j=1,2,\cdots,r\} \label{eq:d9}
\end{eqnarray}
and
\begin{eqnarray}
 T_0(\boldsymbol{x},\boldsymbol{y})=\{\boldsymbol{d}\in R^n\times R^m &|& g_i(\boldsymbol{x},\boldsymbol{y})+\nabla  g_i(\boldsymbol{x},\boldsymbol{y})^\top\boldsymbol{d}\leq 0,i=1,2,\cdots,s,\nonumber\\
&& \nabla h_j(\boldsymbol{x},\boldsymbol{y})^\top\boldsymbol{d}=0, j=1,2,\cdots,r\}.\label{eq:d10}
\end{eqnarray}
It is clear that $T_0(\boldsymbol{x},\boldsymbol{y})\subset T(\boldsymbol{x},\boldsymbol{y})$. For $(\boldsymbol{x}^*,\boldsymbol{y}^*)\in X(g)$, we have
$X(g)\subset [T(\boldsymbol{x}^*,\boldsymbol{y}^*)+(\boldsymbol{x}^*,\boldsymbol{y}^*)].$

Let a  set at a fixed $(\boldsymbol{x}^*,\boldsymbol{y}^*)$ be defined  by
\begin{eqnarray*}
 K_g(\boldsymbol{x}^*,\boldsymbol{y}^*)=\{(\boldsymbol{x},\boldsymbol{y}) \in R^n\times R^m &|& \nabla g(\boldsymbol{x}^*,\boldsymbol{y}^*)^\top [(\boldsymbol{x},\boldsymbol{y})-(\boldsymbol{x}^*,\boldsymbol{y}^*)]<0\}.
\end{eqnarray*}
We have the following Lemmas.

\begin{lemma}
Let $(\boldsymbol{x}^*,\boldsymbol{y}^*)\in X(g)$. If
\begin{eqnarray}
X(f)\cap K_g(\boldsymbol{x}^*,\boldsymbol{y}^*)=\emptyset \label{eq:d11}
\end{eqnarray}
holds, then $g(\boldsymbol{x}^*,\boldsymbol{y}^*)\leq g(\boldsymbol{x},\boldsymbol{y})$ for all $(\boldsymbol{x},\boldsymbol{y})\in X(g)$, i.e. $(\boldsymbol{x}^*,\boldsymbol{y}^*)$  is an optimal solution  to (CNP) and $\boldsymbol{x}^*$ is an optimal solution to (UOP)
\end{lemma}



\begin{lemma}
 Let $(\boldsymbol{x}^*,\boldsymbol{y}^*)\in X(g)$. If
\begin{eqnarray}
X(f)\cap K_g(\boldsymbol{x}^*,\boldsymbol{y}^*)\cap [T(\boldsymbol{x}^*,\boldsymbol{y}^*)\backslash T_0(\boldsymbol{x}^*,\boldsymbol{y}^*)+(\boldsymbol{x}^*,\boldsymbol{y}^*)]=\emptyset \label{eq:d12}
\end{eqnarray}
holds, then $g(\boldsymbol{x}^*,\boldsymbol{y}^*)\leq g(\boldsymbol{x},\boldsymbol{y})$ for all $(\boldsymbol{x},\boldsymbol{y})\in X(g)\cap [T(\boldsymbol{x}^*,\boldsymbol{y}^*)\backslash T_0(\boldsymbol{x}^*,\boldsymbol{y}^*)+(\boldsymbol{x}^*,\boldsymbol{y}^*)]$.
\end{lemma}

\begin{theorem}
 Suppose
\begin{eqnarray}
[(\boldsymbol{x},\boldsymbol{y})-(\boldsymbol{x}^*,\boldsymbol{y}^*)]\in T_0(\boldsymbol{x}^*,\boldsymbol{y}^*), \forall  (\boldsymbol{x},\boldsymbol{y})\in X(g)  \label{eq:d13}
\end{eqnarray}
or \eqref{eq:d12} holds for $(\boldsymbol{x}^*,\boldsymbol{y}^*)\in X(g)$.
If there is an optimal solution $\boldsymbol{d}^*$ to (CNP0)$(\boldsymbol{x}^*,\boldsymbol{y}^*)$ such that $\nabla g(\boldsymbol{x}^*,\boldsymbol{y}^*)^\top \boldsymbol{d}^* \geq 0$, then  $(\boldsymbol{x}^*,\boldsymbol{y}^*)$ is an optimal solution to (CNP) and  $\boldsymbol{x}^*$ is an optimal solution to (UOP) and there are $\boldsymbol{u}^*=(u_1^*,u_2^*,\cdots,u_s^*)^\top\geq 0$ and $\boldsymbol{v}^*=(v_1^*,v_2^*,\cdots,v_r^*)^\top\in R^r$ such that
\begin{eqnarray}
\nabla g(\boldsymbol{x}^*,\boldsymbol{y}^*)+ \nabla \boldsymbol{g}(\boldsymbol{x}^*,\boldsymbol{y}^*)^\top\boldsymbol{u}^* +\nabla \boldsymbol{h}(\boldsymbol{x}^*,\boldsymbol{y}^*)^\top\boldsymbol{v}^*=0.\label{eq:d14}\\
u^*_ig_i(\boldsymbol{x}^*,\boldsymbol{y}^*)^\top=0,i=1,2,\cdots,s.\label{eq:d15}
\end{eqnarray}
\end{theorem}

 {\it Proof.} The proof process is similar to the proof of theorem 7.

\begin{theorem}
Suppose \eqref{eq:d12} or \eqref{eq:d13} holds for $(\boldsymbol{x}^*,\boldsymbol{y}^*)\in X(g)$.
If there are $\boldsymbol{u}^*=(u_1^*,u_2^*,\cdots,u_s^*)^\top\geq 0$ and $\boldsymbol{v}^*=(v_1^*,v_2^*,\cdots,v_r^*)^\top\in R^r$   such that \eqref{eq:d14} and \eqref{eq:d15}  hold, then  $\boldsymbol{x}^*$ is an optimal solution to (UOP).
\end{theorem}

{\it Proof.} Because for $\boldsymbol{u}^*=(u_1^*,u_2^*,\cdots,u_s^*)^\top\geq 0$ and $\boldsymbol{v}^*=(v_1^*,v_2^*,\cdots,v_r^*)^\top\in R^r$  \eqref{eq:d14}) and \eqref{eq:d15}) hold, $\boldsymbol{u}^*$, $\boldsymbol{v}^*$ is a feasible solution to  linear programming (DTCP0)$(\boldsymbol{x}^*,\boldsymbol{y}^*)$. It is clear that $\boldsymbol{d}=0$ is  a feasible solution to  linear programming (CNP0)$(\boldsymbol{x}^*,\boldsymbol{y}^*)$. By the strong dual theorem, there is an optimal solution $\boldsymbol{d}^*$ to (CNP0)$(\boldsymbol{x}^*,\boldsymbol{y}^*)$ such that $\nabla g(\boldsymbol{x}^*,\boldsymbol{y}^*)^\top \boldsymbol{d}^* \geq 0$, \eqref{eq:d14}) and \eqref{eq:d15}) hold. Hence, by Theorem 4, $\boldsymbol{x}^*$ is an optimal solution to (UOP).\\

 The following theorem is concluded by combining Theorem 4 and Theorem 5.

\begin{theorem}
Let $(\boldsymbol{x}^*,\boldsymbol{y}^*)\in X(g)$.
 Then there is an optimal solution $\boldsymbol{d}^*$ to (CNP0)$(\boldsymbol{x}^*,\boldsymbol{y}^*)$ such that $\nabla g(\boldsymbol{x}^*,\boldsymbol{y}^*)^\top \boldsymbol{d}^* \geq 0$ if and only if there are $\boldsymbol{u}^*=(u_1^*,u_2^*,\cdots,u_s^*)^\top\geq 0$ and $\boldsymbol{v}^*=(v_1^*,v_2^*,\cdots,v_r^*)^\top\in R^r$ such that
\eqref{eq:d14}) and \eqref{eq:d15})  hold.
\end{theorem}

Theorem 6 means that $\boldsymbol{x}^*$ doesn't have to be an optimal solution to (UOP) if there is an optimal solution $\boldsymbol{d}^*$ to (CNP0)$(\boldsymbol{x}^*,\boldsymbol{y}^*)$ such that $\nabla g(\boldsymbol{x}^*,\boldsymbol{y}^*)^\top \boldsymbol{d}^* \geq 0$ for $(\boldsymbol{x}^*,\boldsymbol{y}^*)\in X(g)$.

{\bf Example 5} Consider a nonconvex and non-Lipschitz optimization problem:
\begin{eqnarray*}
\mbox{(EX5)} &\min& \ f(x)=|x_1x_2|^\frac{1}{3}+x_1^2+x_2^2\\
  &s.t.& (x_1,x_2)\in R^2.
\end{eqnarray*}
A form of convertible nonconvex optimization of (EX5) is defined by
\begin{eqnarray*}
\mbox{(MEX5-1)} &\min& \ g(x_1,x_2,y_1,y_2,y_3,y_4)=y_4+x_1^2+x_2^2,\\
&s.t.& g_1(x_1,x_2,y_1,y_2,y_3,y_4)=-y_4\leq 0,\\
&&h_1(x_1,x_2,y_1,y_2,y_3,y_4)= 0.5(x_1+x_2)^2-y_1-0.5y_2=0,\\
&&h_2(x_1,x_2,y_1,y_2,y_3,y_4)=x_1^2+x_2^2-y_2=0,\\
&&h_3(x_1,x_2,y_1,y_2,y_3,y_4)= y_1^2-y_3=0,\\
&&h_4(x_1,x_2,y_1,y_2,y_3,y_4)=y_4^6-y_3=0.
\end{eqnarray*}
It is clear that $(\boldsymbol{x}^*,\boldsymbol{y}^*)=(x_1^*,x_2^*,y_1^*,y_2^*,y_3^*,y_4^*)=(0,0,0,0,0,0)\in X(g)$. Let $\boldsymbol{d}=(d_1,d_2,d_3,d_4,d_5,d_6)^\top$ correspond to $(x_1,x_2,y_1,y_2,y_3,y_4)$. We have
\begin{eqnarray*}
\mbox{(MEX5-1)}(0,0,0,0,0,0)\qquad & \min\; & \nabla g(0,0,0,0,0,0)^\top \boldsymbol{d}=d_6\geq 0\\
& \mbox{s.t.}\; &  -d_6\leq 0,-d_3-0.5d_4\leq 0,-d_4\leq 0,-d_5\leq 0,-d_5\leq 0,\\
         &&      d_1,d_2,d_3,d_4,d_5,d_6\in R.
\end{eqnarray*}
By Theorem 1, $(\boldsymbol{x}^*,\boldsymbol{y}^*)=(0,0,0,0,0,0)$ is an optimal solution to (MEX5-1) and
\begin{eqnarray*}
\nabla g(\boldsymbol{x}^*,\boldsymbol{y}^*)+ u_1 \nabla g_1(\boldsymbol{x}^*,\boldsymbol{y}^*)+v_1\nabla h_1(\boldsymbol{x}^*,\boldsymbol{y}^*)+v_2\nabla h_2(\boldsymbol{x}^*,\boldsymbol{y}^*)\\ +v_3\nabla h_3(\boldsymbol{x}^*,\boldsymbol{y}^*)+v_4\nabla h_4(\boldsymbol{x}^*,\boldsymbol{y}^*)=0,
\end{eqnarray*}
where $u_1=1,(v_1,v_2,v_3,v_4)=(0,0,0,0)$.

Example 5  shows that the convertible nonconvex optimization may obtain an optimal solution to (UOP) via the convertible nonconvex forms of $f$ if $f$ is a nonconvex and nonsmooth. If $\min f(\boldsymbol{x})$ is noncovex and nonsmooth optimization, it is necessary to define the subgradient of $f$ in \cite{Huang} such that optimization condition and dual may be obtained. Especially, smoothing function $f$ is necessary if its algorithm is designed. In next section, we discuss that Lagrange dual and algorithm of (CNP) to avoid using smoothing techniques and sub gradients.

\section{Lagrange Dual and Algorithm of (CNP)}

It is well known that Lagrange duality can be used to solve optimization problems.  Especially, there is a zero gap between the optimal objective value of the Lagrangian dual problem of the convex optimization problem and the optimal objective value of the original problem under some constraint qualifications, which is a very important advantage in using the dual problem of convex optimization to get the global optimal solution. In this section, we will establish the Lagrangian duality of (CNP) and its algorithm.

Let $R^s_+=\{\boldsymbol{u}\in R^s| \boldsymbol{u}\geq 0 \}$ and $R^r_+=\{\boldsymbol{v}\in R^r| \boldsymbol{v}\geq 0 \}$. For any $(\boldsymbol{x},\boldsymbol{y})\in R^n\times R^m, \boldsymbol{u}\in R^s_+, \boldsymbol{v}\in R^r$, a Lagrange function of (CNP) is defined by
\begin{eqnarray}
L(\boldsymbol{x},\boldsymbol{y};\boldsymbol{u},\boldsymbol{v})=g(\boldsymbol{x},\boldsymbol{y})+\boldsymbol{u}^\top\boldsymbol{g}(\boldsymbol{x},\boldsymbol{y})+\boldsymbol{v}^\top\boldsymbol{h}(\boldsymbol{x},\boldsymbol{y}). \label{eq:d16}
\end{eqnarray}
Let a dual function of $L(\boldsymbol{x},\boldsymbol{y};\boldsymbol{u},\boldsymbol{v})$ on $(\boldsymbol{u},\boldsymbol{v})\in R^s_+\times R^r$ be defined by
\begin{eqnarray}
\theta(\boldsymbol{u},\boldsymbol{v})= \min \{L(\boldsymbol{x},\boldsymbol{y};\boldsymbol{u},\boldsymbol{v})| (\boldsymbol{x},\boldsymbol{y})\in R^n\times R^m\}. \label{eq:d17}
\end{eqnarray}

For any $(\boldsymbol{x},\boldsymbol{y})\in R^n\times R^m, (\boldsymbol{u},\boldsymbol{v})\in R^s_+\times R^r_+$,  a Lagrange function of (CNP) is defined by
\begin{eqnarray}
L_+(\boldsymbol{x},\boldsymbol{y};\boldsymbol{u},\boldsymbol{v})=g(\boldsymbol{x},\boldsymbol{y})+\boldsymbol{u}^\top\boldsymbol{g}(\boldsymbol{x},\boldsymbol{y})+\boldsymbol{v}^\top\boldsymbol{h}(\boldsymbol{x},\boldsymbol{y}). \label{eq:d18}
\end{eqnarray}
Let a  dual function of $L_+(\boldsymbol{x},\boldsymbol{y};\boldsymbol{u},\boldsymbol{v})$ on $(\boldsymbol{u},\boldsymbol{v})\in R^s_+\times R^r_+$ be defined by
\begin{eqnarray}
\theta_+(\boldsymbol{u},\boldsymbol{v})= \min \{L(\boldsymbol{x},\boldsymbol{y};\boldsymbol{u},\boldsymbol{v})| (\boldsymbol{x},\boldsymbol{y})\in R^n\times R^m\}. \label{eq:d19}
\end{eqnarray}

Let  $(\boldsymbol{x}^*,\boldsymbol{y}^*)\in  R^n\times R^m$ and $(\boldsymbol{u}^*,\boldsymbol{v}^*)\in R^s_+\times R^r$ be given.
For  all $(\boldsymbol{x},\boldsymbol{y})\in R^n\times R^m$ and $(\boldsymbol{u},\boldsymbol{v})\in R^s_+\times R^r$, if
\begin{eqnarray}
L(\boldsymbol{x}^*,\boldsymbol{y}^*;\boldsymbol{u},\boldsymbol{v}) \leq L(\boldsymbol{x}^*,\boldsymbol{y}^*;\boldsymbol{u}^*,\boldsymbol{v}^*)\leq L(\boldsymbol{x},\boldsymbol{y};\boldsymbol{u}^*,\boldsymbol{v}^*), \label{eq:d20}
\end{eqnarray}
then $(\boldsymbol{x}^*,\boldsymbol{y}^*;\boldsymbol{u}^*,\boldsymbol{v}^*)$ is called a saddle point of $L(\boldsymbol{x},\boldsymbol{y};\boldsymbol{u},\boldsymbol{v})$.

Let  $(\boldsymbol{x}^*,\boldsymbol{y}^*)\in  R^n\times R^m$ and $(\boldsymbol{u}^*,\boldsymbol{v}^*)\in R^s_+\times R^r_+$ be given.
For  all $(\boldsymbol{x},\boldsymbol{y})\in R^n\times R^m$ and $(\boldsymbol{u},\boldsymbol{v})\in R^s_+\times R^r_+$, if
\begin{eqnarray}
L_+(\boldsymbol{x}^*,\boldsymbol{y}^*;\boldsymbol{u},\boldsymbol{v}) \leq L_+(\boldsymbol{x}^*,\boldsymbol{y}^*;\boldsymbol{u}^*,\boldsymbol{v}^*)\leq L_+(\boldsymbol{x},\boldsymbol{y};\boldsymbol{u}^*,\boldsymbol{v}^*),\label{eq:d21}
\end{eqnarray}
then $(\boldsymbol{x}^*,\boldsymbol{y}^*;\boldsymbol{u}^*,\boldsymbol{v}^*)$ is called a saddle point of $L_+(\boldsymbol{x},\boldsymbol{y};\boldsymbol{u},\boldsymbol{v})$.

The following conclusions are clear.
\begin{proposition}
If $(\boldsymbol{u},\boldsymbol{v})\in R^s_+\times R^r_+$, then $L(\boldsymbol{x},\boldsymbol{y};\boldsymbol{u},\boldsymbol{v})=L_+(\boldsymbol{x},\boldsymbol{y};\boldsymbol{u},\boldsymbol{v})$ and
$\theta(\boldsymbol{u},\boldsymbol{v})=\theta_+(\boldsymbol{u},\boldsymbol{v})$.
\end{proposition}
\begin{proposition}
If $(\boldsymbol{x}^*,\boldsymbol{y}^*;\boldsymbol{u}^*,\boldsymbol{v}^*)$ is a saddle point of $L(\boldsymbol{x},\boldsymbol{y};\boldsymbol{u},\boldsymbol{v})$ with $R^s_+\times R^r_+$, then $(\boldsymbol{x}^*,\boldsymbol{y}^*;\boldsymbol{u}^*,\boldsymbol{v}^*)$ is  a  saddle point of $L_+(\boldsymbol{x},\boldsymbol{y};\boldsymbol{u},\boldsymbol{v})$.
\end{proposition}

A dual optimization problem of (CNP) is defined by
\begin{eqnarray*}
\mbox{(DCNP)} & \max \ \theta(\boldsymbol{u},\boldsymbol{v}) & \ s.t. \ (\boldsymbol{u},\boldsymbol{v})\in R^s_+\times R^r.
\end{eqnarray*}
It is clear that $\theta(\boldsymbol{u},\boldsymbol{v})$ is a concave function on $(\boldsymbol{u},\boldsymbol{v})\in R^s_+\times R^r$.

A dual optimization problem of (CNP) is defined by
\begin{eqnarray*}
\mbox{(PDCNP)} & \max \ \theta_+(\boldsymbol{u},\boldsymbol{v}) & \ s.t. \ (\boldsymbol{u},\boldsymbol{v})\in R^s_+\times R^r_+.
\end{eqnarray*}
It is clear that $\theta_+(\boldsymbol{u},\boldsymbol{v})$ is a concave function on $R^s_+\times R^r_+$. By \eqref{eq:d16},\eqref{eq:d17},\eqref{eq:d18} and \eqref{eq:d19}, the following weak duality is clear.

\begin{proposition}
 (i) For all $(\boldsymbol{x},\boldsymbol{y})\in R^n\times R^m$ and $(\boldsymbol{u},\boldsymbol{v})\in R^s_+\times R^r$, $L(\boldsymbol{x},\boldsymbol{y};\boldsymbol{u},\boldsymbol{v}) \geq \theta(\boldsymbol{u},\boldsymbol{v})$.

(ii)  For all  $ (\boldsymbol{x},\boldsymbol{y})\in X(g)$ and $(\boldsymbol{u},\boldsymbol{v})\in R^s_+\times R^r$, $g(\boldsymbol{x},\boldsymbol{y})\geq \theta(\boldsymbol{u},\boldsymbol{v})$.

(iii) For all $(\boldsymbol{x},\boldsymbol{y})\in R^n\times R^m$ and $(\boldsymbol{u},\boldsymbol{v})\in R^s_+\times R^r_+$, $L_+(\boldsymbol{x},\boldsymbol{y};\boldsymbol{u},\boldsymbol{v}) \geq \theta_+(\boldsymbol{u},\boldsymbol{v})$.

(iv)  For all $ (\boldsymbol{x},\boldsymbol{y})\in X(g)$ and $(\boldsymbol{u},\boldsymbol{v})\in R^s_+\times R^r_+$,   $g(\boldsymbol{x},\boldsymbol{y})\geq \theta_+(\boldsymbol{u},\boldsymbol{v})$.
\end{proposition}

The strong duality theorem is true as follows.

\begin{theorem}
 (i) If $L(\boldsymbol{x}^*,\boldsymbol{y}^*,\boldsymbol{u}^*)= \theta(\boldsymbol{u}^*,\boldsymbol{v}^*)$ at $(\boldsymbol{x}^*,\boldsymbol{y}^*)\in R^n\times R^m $ and $(\boldsymbol{u}^*,\boldsymbol{v}^*)\in R^s_+\times R^r$, then
\begin{eqnarray}
\nabla L(\boldsymbol{x}^*,\boldsymbol{y}^*;\boldsymbol{u}^*,\boldsymbol{v}^*)=\nabla g(\boldsymbol{x}^*,\boldsymbol{y}^*)+\nabla \boldsymbol{g}(\boldsymbol{x}^*,\boldsymbol{y}^*)^\top \boldsymbol{u}^*+\nabla\boldsymbol{h}(\boldsymbol{x}^*,\boldsymbol{y}^*)^\top\boldsymbol{v}^*=0. \label{eq:d22}
\end{eqnarray}

(ii) If $g(\boldsymbol{x}^*,\boldsymbol{y}^*)= \theta(\boldsymbol{u}^*,\boldsymbol{v}^*)$ at $(\boldsymbol{x}^*,\boldsymbol{y}^*)\in  X(g)$ and $(\boldsymbol{u}^*,\boldsymbol{v}^*)\in R^s_+\times R^r$, then \eqref{eq:d22} and
\begin{eqnarray}
u_i^*g_i(\boldsymbol{x}^*,\boldsymbol{y}^*)=0,i=1,2,\cdots,s \label{eq:d23}
\end{eqnarray}
 hold, $(\boldsymbol{x}^*,\boldsymbol{y}^*)$ is an optimal solution to (CNP) and $(\boldsymbol{u}^*,\boldsymbol{v}^*)$ is an optimal solution to (DCNP).

(iii) Let $(\boldsymbol{x}^*,\boldsymbol{y}^*)\in R^n\times R^m $ and $(\boldsymbol{u}^*,\boldsymbol{v}^*)\in R^s_+\times R^r_+$. Then
\begin{eqnarray}
\nabla L_+(\boldsymbol{x}^*,\boldsymbol{y}^*;\boldsymbol{u}^*,\boldsymbol{v}^*)=\nabla g(\boldsymbol{x}^*,\boldsymbol{y}^*)+\nabla \boldsymbol{g}(\boldsymbol{x}^*,\boldsymbol{y}^*)^\top \boldsymbol{u}^*+\nabla\boldsymbol{h}(\boldsymbol{x}^*,\boldsymbol{y}^*)^\top\boldsymbol{v}^*=0 \label{eq:d24}
\end{eqnarray}
holds if and only if $L_+(\boldsymbol{x}^*,\boldsymbol{y}^*,\boldsymbol{u}^*,\boldsymbol{v}^*)= \theta_+(\boldsymbol{u}^*,\boldsymbol{v}^*)$.

(iv) Let $(\boldsymbol{x}^*,\boldsymbol{y}^*)\in  X(f)$ and $(\boldsymbol{u}^*,\boldsymbol{v}^*)\in R^s_+\times R^r_+$. Then \eqref{eq:d24} and
\begin{eqnarray}
u_i^*g_i(\boldsymbol{x}^*,\boldsymbol{y}^*)=0,i=1,2,\cdots,s \label{eq:d25}
\end{eqnarray}
 hold if and only if $g(\boldsymbol{x}^*,\boldsymbol{y}^*)= \theta_+(\boldsymbol{u}^*,\boldsymbol{v}^*)$. Furthermore, $(\boldsymbol{x}^*,\boldsymbol{y}^*)$ is an optimal solution to (CNP) and $(\boldsymbol{u}^*,\boldsymbol{v}^*)$ is an optimal solution to (PDCNP).
\end{theorem}

By Theorem 7, there is a zero gap between the optimal objective value of the Lagrangian dual problem of (CNP) and the optimal objective value of the original problem (CNP). But, if a Lagrange function $L(\boldsymbol{x},\boldsymbol{y};\boldsymbol{u}^*,\boldsymbol{v}^*)$ is convex on $(\boldsymbol{x},\boldsymbol{y})$  at $(\boldsymbol{u}^*,\boldsymbol{v}^*)\in R^s_+\times R^r$, the conclusion of zero gap holds by the proof of Theorem 7(iii) and (iv).

\begin{corollary}
(i) Suppose that for $(\boldsymbol{x}^*,\boldsymbol{y}^*)\in R^n\times R^m $ and $(\boldsymbol{u}^*,\boldsymbol{v}^*)\in R^s_+\times R^r$, the Lagrange function $L(\boldsymbol{x},\boldsymbol{y};\boldsymbol{u}^*,\boldsymbol{v}^*)$ is convex on $(\boldsymbol{x},\boldsymbol{y})$  at $(\boldsymbol{u}^*,\boldsymbol{v}^*)$. If \eqref{eq:d23} holds, then $L(\boldsymbol{x}^*,\boldsymbol{y}^*,$ $\boldsymbol{u}^*)= \theta(\boldsymbol{u}^*,\boldsymbol{v}^*)$.

(ii) Suppose that for $(\boldsymbol{x}^*,\boldsymbol{y}^*)\in  X(f)$ and $(\boldsymbol{u}^*,\boldsymbol{v}^*)\in R^s_+\times R^r$,  the Lagrange function $L(\boldsymbol{x},\boldsymbol{y};\boldsymbol{u}^*,\boldsymbol{v}^*)$ is convex on $(\boldsymbol{x},\boldsymbol{y})$  at $(\boldsymbol{u}^*,\boldsymbol{v}^*)$. If \eqref{eq:d23}   holds,
 then $g(\boldsymbol{x}^*,\boldsymbol{y}^*)= \theta(\boldsymbol{u}^*,\boldsymbol{v}^*)$.
\end{corollary}

\begin{proposition}
 (i) Suppose that $(\boldsymbol{x}^*,\boldsymbol{y}^*)\in R^n\times R^m $ and $(\boldsymbol{u}^*,\boldsymbol{v}^*)\in R^s_+\times R^r$. Then
$(\boldsymbol{x}^*,\boldsymbol{y}^*,\boldsymbol{u}^*,\boldsymbol{v}^*)$ is a saddle point of $L(\boldsymbol{x},\boldsymbol{y},\boldsymbol{u},\boldsymbol{v})$
if and only if $L(\boldsymbol{x}^*,\boldsymbol{y}^*,\boldsymbol{u}^*,\boldsymbol{v}^*)= \theta(\boldsymbol{u}^*,\boldsymbol{v}^*)$.

(ii) Suppose that \eqref{eq:d23} holds at $(\boldsymbol{x}^*,\boldsymbol{y}^*)\in  X(f)$ and $(\boldsymbol{u}^*,\boldsymbol{v}^*)\in R^s_+\times R^r$. Then $(\boldsymbol{x}^*,\boldsymbol{y}^*;\boldsymbol{u}^*,\boldsymbol{v}^*)$ is a saddle point of $L(\boldsymbol{x},\boldsymbol{y},\boldsymbol{u},\boldsymbol{v})$ if and only if $g(\boldsymbol{x}^*,\boldsymbol{y}^*)= \theta(\boldsymbol{u}^*,\boldsymbol{v}^*)$.

(iii) Suppose that $(\boldsymbol{x}^*,\boldsymbol{y}^*)\in R^n\times R^m $ and $(\boldsymbol{u}^*,\boldsymbol{v}^*)\in R^s_+\times R^r_+$. Then
$(\boldsymbol{x}^*,\boldsymbol{y}^*,\boldsymbol{u}^*,\boldsymbol{v}^*)$ is a saddle point of $L_+(\boldsymbol{x},\boldsymbol{y},\boldsymbol{u},\boldsymbol{v})$
if and only if $L_+(\boldsymbol{x}^*,\boldsymbol{y}^*,\boldsymbol{u}^*,\boldsymbol{v}^*)= \theta_+(\boldsymbol{u}^*,\boldsymbol{v}^*)$.

(iv) Suppose that \eqref{eq:d25} holds at $(\boldsymbol{x}^*,\boldsymbol{y}^*)\in  X(f)$ and $(\boldsymbol{u}^*,\boldsymbol{v}^*)\in R^s_+\times R^r_+$. Then $(\boldsymbol{x}^*,\boldsymbol{y}^*;\boldsymbol{u}^*,\boldsymbol{v}^*)$ is a saddle point of $L_+(\boldsymbol{x},\boldsymbol{y},\boldsymbol{u},\boldsymbol{v})$ if and only if $g(\boldsymbol{x}^*,\boldsymbol{y}^*)= \theta_+(\boldsymbol{u}^*,\boldsymbol{v}^*)$.\\
\end{proposition}

Define an augmented Lagrange penalty function of (CNP) as follows
\begin{eqnarray}
A(\boldsymbol{x},\boldsymbol{y};\boldsymbol{u},\boldsymbol{v},\rho)&=&g(\boldsymbol{x},\boldsymbol{y})+
\boldsymbol{u}^\top\boldsymbol{g}(\boldsymbol{x},\boldsymbol{y})+
\boldsymbol{v}^\top\boldsymbol{h}(\boldsymbol{x},\boldsymbol{y})\nonumber\\ && +\rho \sum\limits_{i=1}^s g_i^+(\boldsymbol{x},\boldsymbol{y})^2
+\rho \sum\limits_{j=1}^r h_j(\boldsymbol{x},\boldsymbol{y})^2, \nonumber\\ &&
 (\boldsymbol{x},\boldsymbol{y})\in R^n\times R^m,(\boldsymbol{u},\boldsymbol{v})\in R^s_+\times R^r, \label{eq:d26}
\end{eqnarray}
where $\rho>0$ is a penalty parameter and $g_i^+(\boldsymbol{x},\boldsymbol{y})=\max\{g_i(\boldsymbol{x},\boldsymbol{y}),0\}$.
Define an  augmented Lagrange penalty function of (CNP) as follows
\begin{eqnarray}
A_+(\boldsymbol{x},\boldsymbol{y};\boldsymbol{u},\boldsymbol{v},\rho)&=&g(\boldsymbol{x},\boldsymbol{y})
+\boldsymbol{u}^\top\boldsymbol{g}(\boldsymbol{x},\boldsymbol{y})
+\boldsymbol{v}^\top\boldsymbol{h}(\boldsymbol{x},\boldsymbol{y})\nonumber\\&& +\rho \sum\limits_{i=1}^s g_i^+(\boldsymbol{x},\boldsymbol{y})^2 +\rho \sum\limits_{j=1}^r h_j(\boldsymbol{x},\boldsymbol{y})^2,\nonumber\\
&& \ \ (\boldsymbol{x},\boldsymbol{y})\in R^n\times R^m, (\boldsymbol{u},\boldsymbol{v})\in R^s_+\times R^r_+, \label{eq:d27}
\end{eqnarray}
where $\rho>0$ is a penalty parameter and $g_i^+(\boldsymbol{x},\boldsymbol{y})=\max\{g_i(\boldsymbol{x},\boldsymbol{y}),0\}$. Theorem 8 shows that \eqref{eq:d26} and \eqref{eq:d27} are exact penalty function for $\rho>0$.

\begin{theorem}
 (i) If $(\boldsymbol{x}^*,\boldsymbol{y}^*)\in X(f)$ is an optimal solution to $\min\limits_{(\boldsymbol{x},\boldsymbol{y})}L(\boldsymbol{x},\boldsymbol{y};\boldsymbol{u}^*,\boldsymbol{v}^*)$ at $(\boldsymbol{u}^*,\boldsymbol{v}^*)\in R^s_+\times R^r$,  then $(\boldsymbol{x}^*,\boldsymbol{y}^*)$ is an optimal solution to $\min\limits_{(\boldsymbol{x},\boldsymbol{y})}A(\boldsymbol{x},\boldsymbol{y};\boldsymbol{u}^*,\boldsymbol{v}^*, \rho)$ for all $\rho>0$.

(ii) Suppose that \eqref{eq:d25} holds for $(\boldsymbol{x}^*,\boldsymbol{y}^*)\in X(f)$ and $(\boldsymbol{u}^*,\boldsymbol{v}^*)\in R^s_+\times R^r_+$. Then \eqref{eq:d24} holds
  if and only if $(\boldsymbol{x}^*,\boldsymbol{y}^*)$ is an optimal solution to $\min\limits_{(\boldsymbol{x},\boldsymbol{y})}A_+(\boldsymbol{x},\boldsymbol{y};\boldsymbol{u}^*,\boldsymbol{v}^*, \rho)$ for all $\rho>0$.
\end{theorem}

{\it Proof.}  (i) For any $(\boldsymbol{x},\boldsymbol{y})\in R^n\times R^m$, we have
\begin{eqnarray*}
A(\boldsymbol{x},\boldsymbol{y};\boldsymbol{u}^*,\boldsymbol{v}^*,\rho)&=&g(\boldsymbol{x},\boldsymbol{y})+\sum\limits_{i=1}^s u_i^*g_i(\boldsymbol{x},\boldsymbol{y})+\rho \sum\limits_{i=1}^s g_i^+(\boldsymbol{x},\boldsymbol{y})^2\\&& +\sum\limits_{i=1}^r v_j^*h_j(\boldsymbol{x},\boldsymbol{y})+\rho \sum\limits_{j=1}^r h_j(\boldsymbol{x},\boldsymbol{y})^2 \\
&\geq &L(\boldsymbol{x},\boldsymbol{y};\boldsymbol{u}^*,\boldsymbol{v}^*)+\rho \sum\limits_{i=1}^s g_i^+(\boldsymbol{x},\boldsymbol{y})^2+\rho \sum\limits_{j=1}^r h_j(\boldsymbol{x},\boldsymbol{y})^2 \\
&\geq &L(\boldsymbol{x}^*,\boldsymbol{y}^*;\boldsymbol{u}^*,\boldsymbol{v}^*)=A(\boldsymbol{x}^*,\boldsymbol{y}^*;\boldsymbol{u}^*,\boldsymbol{v}^*,\rho).
\end{eqnarray*}
Hence, $(\boldsymbol{x}^*,\boldsymbol{y}^*)$ is an optimal solution to $\min\limits_{(\boldsymbol{x},\boldsymbol{y})}A(\boldsymbol{x},\boldsymbol{y};\boldsymbol{u}^*,\boldsymbol{v}^*, \rho)$ for all $\rho>0$.

(ii) Let us prove that $(\boldsymbol{x}^*,\boldsymbol{y}^*)$ is an optimal solution to $\min\limits_{(\boldsymbol{x},\boldsymbol{y})}A_+(\boldsymbol{x},\boldsymbol{y};\boldsymbol{u}^*,\boldsymbol{v}^*, \rho)$ for all $\rho>0$, when  \eqref{eq:d24} and \eqref{eq:d25} hold. Let any $(\boldsymbol{x},\boldsymbol{y})\in R^n\times R^m$. Since $g(\boldsymbol{x},\boldsymbol{y})$, $g_i(\boldsymbol{x},\boldsymbol{y})(i=1,2,\cdots,s)$ and $h_j(\boldsymbol{x},\boldsymbol{y})(i=1,2,\cdots,r)$ are convex, we have
\begin{eqnarray}
g(\boldsymbol{x},\boldsymbol{y})-g(\boldsymbol{x}^*,\boldsymbol{y}^*)\geq \nabla g(\boldsymbol{x}^*,\boldsymbol{y}^*)[(\boldsymbol{x},\boldsymbol{y})- (x^*,\boldsymbol{y}^*)],\label{eq:d28}\\
g_i(\boldsymbol{x},\boldsymbol{y})-g_i(\boldsymbol{x}^*,\boldsymbol{y}^*)\geq \nabla g_i(\boldsymbol{x}^*,\boldsymbol{y}^*)[(\boldsymbol{x},\boldsymbol{y})- (\boldsymbol{x}^*,\boldsymbol{y}^*)],i=1,2,\cdots,s,\label{eq:d29}\\
h_j(\boldsymbol{x},\boldsymbol{y})-h_j(\boldsymbol{x}^*,\boldsymbol{y}^*)\geq \nabla h_j(\boldsymbol{x}^*,\boldsymbol{y}^*)[(\boldsymbol{x},\boldsymbol{y})- (\boldsymbol{x}^*,\boldsymbol{y}^*)],i=1,2,\cdots,r.\label{eq:d30}
\end{eqnarray}
From \eqref{eq:d28} , \eqref{eq:d29}, \eqref{eq:d30} and Theorem7(iv), we have
\begin{eqnarray*}
A_+(\boldsymbol{x},\boldsymbol{y};\boldsymbol{u}^*,\boldsymbol{v}^*, \rho)&=&g(\boldsymbol{x},\boldsymbol{y})+\sum\limits_{i=1}^s u_i^*g_i(\boldsymbol{x},\boldsymbol{y})+\rho \sum\limits_{i=1}^s g_i^+(\boldsymbol{x},\boldsymbol{y})^2\\ && +\sum\limits_{i=1}^r v_j^*h_j(\boldsymbol{x},\boldsymbol{y})+\rho \sum\limits_{j=1}^r h_j(\boldsymbol{x},\boldsymbol{y})^2 \\
&\geq &L_+(\boldsymbol{x},\boldsymbol{y};\boldsymbol{u}^*,\boldsymbol{v}^*)+\rho \sum\limits_{i=1}^s g_i^+(\boldsymbol{x},\boldsymbol{y})^2+\rho \sum\limits_{j=1}^r h_j(\boldsymbol{x},\boldsymbol{y})^2 \\
&\geq &g(\boldsymbol{x}^*,\boldsymbol{y}^*)=A_+(\boldsymbol{x}^*,\boldsymbol{y}^*;\boldsymbol{u}^*,\boldsymbol{v}^*,\rho).
\end{eqnarray*}
Hence, $(\boldsymbol{x}^*,\boldsymbol{y}^*)$ is an optimal solution to $\min\limits_{(\boldsymbol{x},\boldsymbol{y})}A_+(\boldsymbol{x},\boldsymbol{y};\boldsymbol{u}^*,\boldsymbol{v}^*, \rho)$ for all $\rho>0$.

Now, when $(\boldsymbol{x}^*,\boldsymbol{y}^*)$ is an optimal solution to $\min\limits_{(\boldsymbol{x},\boldsymbol{y})}A_+(\boldsymbol{x},\boldsymbol{y};\boldsymbol{u}^*,\boldsymbol{v}^*, \rho)$ for all $\rho>0$, it is clear that
\begin{eqnarray*}
\nabla A_+(\boldsymbol{x}^*,\boldsymbol{y}^*;\boldsymbol{u}^*,\boldsymbol{v}^*,\rho)&=&\nabla g(\boldsymbol{x}^*,\boldsymbol{y}^*)+\sum\limits_{i=1}^s u_i^*\nabla g_i(\boldsymbol{x}^*,\boldsymbol{y}^*)\\ && +\rho \sum\limits_{i=1}^s 2g_i^+(\boldsymbol{x}^*,\boldsymbol{y}^*) \nabla g_i(\boldsymbol{x}^*,\boldsymbol{y}^*)
 +\sum\limits_{j=1}^r v_i^*\nabla h_j(\boldsymbol{x}^*,\boldsymbol{y}^*)\\ && +\rho \sum\limits_{j=1}^r 2h_j(\boldsymbol{x}^*,\boldsymbol{y}^*) \nabla h_j(\boldsymbol{x}^*,\boldsymbol{y}^*), \\
&=&\nabla L_+(\boldsymbol{x}^*,\boldsymbol{y}^*;\boldsymbol{u}^*,\boldsymbol{v}^*)=0.
\end{eqnarray*}
Hence, $L_+(\boldsymbol{x}^*,\boldsymbol{y}^*;\boldsymbol{u}^*,\boldsymbol{v}^*)=\theta(\boldsymbol{u}^*,\boldsymbol{v}^*)$.\\

Theorem 7 and 8 mean that $g(\boldsymbol{x}^*,\boldsymbol{y}^*)= \theta(\boldsymbol{u}^*,\boldsymbol{v}^*)$ or $g(\boldsymbol{x}^*,\boldsymbol{y}^*)= \theta_+(\boldsymbol{u}^*,\boldsymbol{v}^*)$ does not hold, if there is no  $(\boldsymbol{x}^*,\boldsymbol{y}^*)\in R^n\times R^m $ and $(\boldsymbol{u}^*,\boldsymbol{v}^*)\in R^s_+\times R^r_+$ such that \eqref{eq:d23} or \eqref{eq:d25} holds. This means that  \eqref{eq:d23} or \eqref{eq:d25}  is necessary, if an global optimal solution to (CNP) is to be found by Theorem 7. Let us see the following example.

{\bf Example 6} Consider the optimization(Example 5):
\begin{eqnarray*}
\mbox{(EX5)} &\min& \ f(x)=|x_1x_2|^\frac{1}{3}+x_1^2+x_2^2\\
  &s.t.& (x_1,x_2)\in R^2.
\end{eqnarray*}
The convertible nonconvex optimization of (EX5) is defined by (MEX5-1). The  Lagrange function of (MEX5-1) is defined as
\begin{eqnarray*}
L_+(\boldsymbol{x},\boldsymbol{y};\boldsymbol{u},\boldsymbol{v})=y_4+x_1^2+x_2^2-u_1y_4+v_1(0.5(x_1+x_2)^2-y_1-0.5y_2)
\\ +v_2(x_1^2+x_2^2-y_2)+v_3(y_1^2-y_3)+v_4(y_4^6-y_3),
\end{eqnarray*}
where $\forall u_1\geq 0$ and $(\boldsymbol{x},\boldsymbol{y};\boldsymbol{u},\boldsymbol{v})\in R^2\times R^4\times R^1_+\times R^4_+$.
The  dual function of $L_+(\boldsymbol{x},\boldsymbol{y};\boldsymbol{u},\boldsymbol{v})$ is obtained by
$$
\theta_+(\boldsymbol{u},\boldsymbol{v})=\left\{
\begin{array}{ll}
0 & \mbox{if }~  u_1=1,v_1=v_2=v_3=v_4=0\\
-\infty                       & \mbox{if }~  u_1\not=1\ \mbox{or}\ \boldsymbol{v}\in R^4_+\setminus \{0\}.
\end{array}\right.$$
It is clear that $g(\boldsymbol{x},\boldsymbol{y})\geq \theta_+(\boldsymbol{u},\boldsymbol{v})$ for all $$(\boldsymbol{x},\boldsymbol{y};\boldsymbol{u},\boldsymbol{v})\in R^2\times R^4\times R^1_+\times R^4_+.$$ When $u^*_1=1$, $v_1^*=v_2^*=v_3^*=v_4^*=0$ and $(\boldsymbol{x}^*,\boldsymbol{y}^*)=(0,0,0,0,0,0)$,  $g(\boldsymbol{x}^*, \boldsymbol{y}^*)= \theta(\boldsymbol{u}^*,\boldsymbol{v}^*)=0$. $(\boldsymbol{x}^*,\boldsymbol{y}^*)$ is the optimal solution to (MEX5-1).
The Lagrange function of (MEX5-1) is defined as
\begin{eqnarray*}
L(\boldsymbol{x},\boldsymbol{y};\boldsymbol{u},\boldsymbol{v})=y_4+x_1^2+x_2^2-u_1y_4+v_1(0.5(x_1+x_2)^2-y_1-0.5y_2)\\
+v_2(x_1^2+x_2^2-y_2)+v_3(y_1^2-y_3)+v_4(y_4^6-y_3),
\end{eqnarray*}
where  $(\boldsymbol{x},\boldsymbol{y};\boldsymbol{u},\boldsymbol{v})\in R^2\times R^4\times R^1_+\times R^4$.
The dual function of $L(\boldsymbol{x},\boldsymbol{y};\boldsymbol{u},\boldsymbol{v})$ is obtained by
$$
\theta_+(\boldsymbol{u},\boldsymbol{v})=\left\{
\begin{array}{ll}
0 & \mbox{if }~  u_1=1,v_1=v_2=v_3=v_4=0\\
-\infty                       & \mbox{if }~  u_1\not=1\ \mbox{or}\ \boldsymbol{v}\in R^4\setminus \{0\}.
\end{array}\right.$$

The  augmented Lagrange penalty function of (MEX5-1) is defined as
\begin{eqnarray*}
A_+(\boldsymbol{x},\boldsymbol{y};\boldsymbol{u},\boldsymbol{v},\rho)&=&y_4+x_1^2+x_2^2-u_1y_4+v_1(0.5(x_1+x_2)^2-y_1-0.5y_2)
\\ && +v_2(x_1^2+x_2^2-y_2)+v_3(y_1^2-y_3)+v_4(y_4^6-y_3)\\&&  +\rho[\max\{-y_4,0\}^2+(0.5(x_1+x_2)^2-y_1-0.5y_2)^2\\
&& +(x_1^2+x_2^2-y_2)^2+(y_1^2-y_3)^2+(y_4^6-y_3)^2],
\end{eqnarray*}
where $\rho>0$ and $(\boldsymbol{x},\boldsymbol{y};\boldsymbol{u},\boldsymbol{v})\in R^2\times R^4\times R^1_+\times R^4_+$. By Theorem 8 and Example 5, it is clear that  $(\boldsymbol{x}^*,\boldsymbol{y}^*)=(0,0,0,0,0,0)$ is an optimal solution to $\min\limits_{(\boldsymbol{x},\boldsymbol{y})}A_+(\boldsymbol{x},\boldsymbol{y};\boldsymbol{u}^*,\boldsymbol{v}^*, \rho)$ for all $\rho>0$ and $(\boldsymbol{u}^*,\boldsymbol{v}^*)=(1,0,0,0,0)$.\\

In order to obtain a solution to (UOP), by Theorem 8, we may  find an approximate solution to (CNP) by the augmented Lagrange penalty optimization as follows
\begin{eqnarray*}
\mbox{(CNP)}(\boldsymbol{u},\boldsymbol{v},\rho) &\min& A(\boldsymbol{x},\boldsymbol{y};\boldsymbol{u},\boldsymbol{v},\rho), \\
                        &s.t.&(\boldsymbol{x},\boldsymbol{y})\in R^n\times R^m.
\end{eqnarray*}

Let $\boldsymbol{g}^+(\boldsymbol{x},\boldsymbol{y})=(g^+_1(\boldsymbol{x},\boldsymbol{y}),g^+_2(\boldsymbol{x},\boldsymbol{y}),\cdots,g^+_s(\boldsymbol{x},\boldsymbol{y}))^\top$.

To solve  (CNP)$(\boldsymbol{u},\boldsymbol{v},\rho)$, we propose an algorithm of augmented Lagrange penalty function of (CNP)(which is called ALPF Algorithm).

{\bf ALPF Algorithm:}

\begin{description}
\item[Step 1:] Let $\epsilon>0,\rho_{1}>0,N>1, (\boldsymbol{x}^0,\boldsymbol{y}^0)\in R^n\times R^m$, $\boldsymbol{v}^1\in R^r$, $\boldsymbol{u}^1=\boldsymbol{g}^+(\boldsymbol{x}^0,\boldsymbol{y}^0)$,$\boldsymbol{v}^{1}=0$, and $k=1$.

\item[Step 2:] Solve $(\boldsymbol{x}^k,\boldsymbol{y}^k)\in R^n\times R^m $ to problem $\lim\limits_{(\boldsymbol{x},\boldsymbol{y})}A(\boldsymbol{x},\boldsymbol{y},\boldsymbol{u}^{k},\boldsymbol{v}^{k},\rho_{k})$ such that  $\nabla A(\boldsymbol{x}^{k},\boldsymbol{y}^{k},\boldsymbol{u}^{k},\boldsymbol{v}^{k},\rho_{k})=0$, and go to Step 3.

\item[Step 3:]If $u_i^kg_i(\boldsymbol{x}^k,\boldsymbol{y}^k)=0(i=1,2,\cdots,s)$ at $(\boldsymbol{x}^k,\boldsymbol{y}^k)\in X(f)$ and $L(\boldsymbol{x},\boldsymbol{y},\boldsymbol{u}^{k},\boldsymbol{v}^k)$ is convex on $(\boldsymbol{x},\boldsymbol{y})$ for $(\boldsymbol{u}^k,\boldsymbol{v}^k)\in R^s_+\times R^r$,  then stop and  $\boldsymbol{x}^k$ is an optimal solution to (UOP). Otherwise, go to Step 4.

\item[Step 4:] If $|A(\boldsymbol{x}^k,\boldsymbol{y}^k,\boldsymbol{u}^{k},\boldsymbol{v}^k,\rho_{k})-g(\boldsymbol{x}^k,\boldsymbol{y}^k)| <\epsilon$ and $\|\boldsymbol{g}^+(\boldsymbol{x}^k,\boldsymbol{y}^k)\|+\|\boldsymbol{h}(\boldsymbol{x}^k,\boldsymbol{y}^k)\|<\epsilon$, then stop and  $\boldsymbol{x}^k$ is an approximate solution to (UOP). Otherwise, for $i=1,2,\cdots,s,$
     let \\
     \kgg{0.3in}$ u^{k+1}_i=$\begin{math} \left\{
\begin{array}{l}
u^{k}_i+2\rho_{k}g^+_i(\boldsymbol{x}^{k},\boldsymbol{y}^{k}), \kgg{0.5in} if\kgg{0.3in} g_i(\boldsymbol{x}^{k},\boldsymbol{y}^{k})\geq 0,\\
0,\kgg{1.7in} if\kg g_i(\boldsymbol{x}^{k},\boldsymbol{y}^{k})< 0,\\
\end{array}\right.\end{math}\\
     $\boldsymbol{v}^{k+1}=\boldsymbol{v}^{k}+2\rho_{k}\boldsymbol{h}(\boldsymbol{x}^{k},\boldsymbol{y}^{k})$, $\rho_{k+1}=N\rho_{k}$, $k:=k+1$ and go to Step 2.
\end{description}

{Note:} By Theorem 7 and Theorem 8, if $A(\boldsymbol{x}^k,\boldsymbol{y}^k;\boldsymbol{u}^{k},\boldsymbol{v}^{k},\rho_{k})=g(\boldsymbol{x}^k,\boldsymbol{y}^k)$ and $(\boldsymbol{x}^k,\boldsymbol{y}^k)\in X(g)$ for some $k$, then $\theta(\boldsymbol{u}^{k},\boldsymbol{v}^{k})=g(\boldsymbol{x}^k,\boldsymbol{y}^k)$ and $\boldsymbol{x}^k$ is an optimal solution to (UOP).
Hence, $\boldsymbol{x}^k$ may be an approximate solution to (UOP) if $|A(\boldsymbol{x}^k,\boldsymbol{y}^k;\boldsymbol{u}^{k},\boldsymbol{v}^{k},\rho_{k})-g(\boldsymbol{x}^k,\boldsymbol{y}^k)| <\epsilon$ and $\|\boldsymbol{g}^+(\boldsymbol{x}^k,\boldsymbol{y}^k)\|+\|\boldsymbol{h}(\boldsymbol{x}^k,\boldsymbol{y}^k)\|<\epsilon$ hold. ALPF Algorithm may be able to find an approximate global optimal solution to (UOP). Under some conditions, that ALPF Algorithm  can converge to a KKT point of (CNP) for $\epsilon=0$ is proved.

Let
$$S(\pi, g)=\{(\boldsymbol{x},\boldsymbol{y})\mid \pi\geq  g(\boldsymbol{x},\boldsymbol{y})\}, $$
which is called a level set.  If $S(\pi, g)$ is bounded for any given $\pi>0$, then $S(\pi, g)$ is called to be bounded.

\begin{lemma}
Suppose that sequence $\{(\boldsymbol{x}^k,\boldsymbol{y}^k)\}$ is obtained by  ALPF Algorithm for every
$k=1,2,\cdots$. Let sequence $\{(\boldsymbol{x}^k,\boldsymbol{y}^k)\}$ converge to $(\boldsymbol{x}^*,\boldsymbol{y}^*)$. If $g_i(\boldsymbol{x}^*,\boldsymbol{y}^*)<0$ for some $i$, then $u_i^k\to 0$ as  $k\to +\infty$.\\
\end{lemma}

{\it Proof.} If $g_i(\boldsymbol{x}^*,\boldsymbol{y}^*)<0$ for some $i$, there is $k'>0$ such that $g_i^k(\boldsymbol{x}^*,\boldsymbol{y}^*)<0$ for all $k>k'$. By Step 4 of  ALPF Algorithm, we have $ u^{k+1}_i=0$ for all $k>k'$. Hence,  $u_i^k\to 0$ as  $k\to +\infty$.

\begin{theorem}
Let $\epsilon=0$. Suppose that sequence $\{(\boldsymbol{x}^k,\boldsymbol{y}^k)\}$ is obtained by  ALPF Algorithm  for every
$k=1,2,\cdots$.
Let, for every $k=1,2,\cdots,$ sequence $\{H_k(\boldsymbol{x}^{k},\boldsymbol{y}^{k},\rho_{k})\}$ be bounded and the level set $S(\pi, g)$ be bounded, where
$$H_k(\boldsymbol{x}^{k},\boldsymbol{y}^{k},\rho_{k})=g(\boldsymbol{x}^{k},\boldsymbol{y}^{k})+\rho_{k}\sum\limits_{i=1}^s g_i^+(\boldsymbol{x}^{k},\boldsymbol{y}^{k})^2+\rho_{k}\sum\limits_{j=1}^r h_j(\boldsymbol{x}^{k},\boldsymbol{y}^{k})^2.$$

(i) If the algorithm stops at finite step $k$, then $\boldsymbol{x}^k$  is a global optimal solution to (UOP).

(ii) If sequence $\{(\boldsymbol{x}^{k},\boldsymbol{y}^{k})\}$ is infinite, then $\{(\boldsymbol{x}^{k},\boldsymbol{y}^{k})\}$ is bounded and any limit point $(\boldsymbol{x}^*,{\boldsymbol{y}}^*)$ of it belongs to $X(f)$, and
there exist $\eta>0$, $\alpha_i\geq 0(i=1, 2, \cdots, s)$ and $\beta_j(j=1, 2, \cdots, r)$, such that
\begin{eqnarray}
\eta\nabla g(\boldsymbol{x}^*,{\boldsymbol{y}}^*)+ \sum\limits_{i=1}^s \alpha_i \nabla g_i(\boldsymbol{x}^*,{\boldsymbol{y}}^*)+\sum\limits_{i=1}^r \beta_j \nabla h_j(\boldsymbol{x}^*,{\boldsymbol{y}}^*)=0. \label{eq:d31}\\
\alpha_i \nabla g_i(\boldsymbol{x}^*,{\boldsymbol{y}}^*)=0,i=1,2,\cdots,s. \label{eq:d32}
\end{eqnarray}
 If $\beta_j\geq 0(j=1, 2, \cdots, r)$,  then
$\boldsymbol{x}^*$ is an optimal solution to (UOP)\\
\end{theorem}

{\it Proof.} (i) The conclusion is clear by Theorem 7 and Corollary 3.

(ii) By ALPF Algorithm,  since $\{H_k(\boldsymbol{x}^k,{\boldsymbol{y}}^k,\rho_k)\}$ is bounded as $k\to +\infty$,  there must be some $\pi>0$ such
that\begin{eqnarray*}
          \pi&>&H_k(\boldsymbol{x}^k,{\boldsymbol{y}}^k,\rho_k)\\
          &=&  g(\boldsymbol{x}^k,{\boldsymbol{y}}^k)+\rho_k\sum\limits_{i=1}^s g_i^+(\boldsymbol{x}^k,{\boldsymbol{y}}^k)^2+\rho_k\sum\limits_{j=1}^r h_j(\boldsymbol{x}^k,{\boldsymbol{y}}^k)^2\\
          &\geq&  g(\boldsymbol{x}^k,{\boldsymbol{y}}^k).
\end{eqnarray*}
 $\{(\boldsymbol{x}^k,{\boldsymbol{y}}^k)\}$ is bounded because the level set $S(\pi,f)$ is bounded.
Without loss of generality, suppose $(\boldsymbol{x}^k,{\boldsymbol{y}}^k)\to (\boldsymbol{x}^*,{\boldsymbol{y}}^*)$.  Since $g$ is continuous, $S(\pi,f)$ is closed. So, $g(\boldsymbol{x}^k,{\boldsymbol{y}}^k)$ is bounded and there is a $\sigma>0$ such that $g(\boldsymbol{x}^k,{\boldsymbol{y}}^k> -\sigma$.

From the above inequality, we have
that\begin{eqnarray*}
       \sum\limits_{i=1}^s g_i^+(\boldsymbol{x}^k,{\boldsymbol{y}}^k)^2+\sum\limits_{j=1}^r h_j(\boldsymbol{x}^k,{\boldsymbol{y}}^k)^2 \leq\frac{1}{\rho_k} (\pi- g(\boldsymbol{x}^k,{\boldsymbol{y}}^k))<\frac{\pi+\sigma}{\rho_k}.
\end{eqnarray*}
 And $\sum\limits_{i=1}^s g_i^+(\boldsymbol{x}^k,{\boldsymbol{y}}^k)^2+\sum\limits_{j=1}^r (h_j(\boldsymbol{x}^k,{\boldsymbol{y}}^k))^2\to 0$ as  $\rho_k\to +\infty$.
So, $(\boldsymbol{x}^*,{\boldsymbol{y}}^*)\in X(f)$.

By  ALPF algorithm,  there is infinite sequence $\{\boldsymbol{x}^k,{\boldsymbol{y}}^k,\boldsymbol{u}^k,\boldsymbol{v}^k,\rho_k) \}$ such that $ \nabla A(\boldsymbol{x}^k,{\boldsymbol{y}}^k;\boldsymbol{u}^k,\boldsymbol{v}^k,\rho_k)=0$. We have
 \begin{eqnarray}
\nabla g(\boldsymbol{x}^k,{\boldsymbol{y}}^k)
+ \sum\limits_{i=1}^s \bar{u}^{k}_i \nabla g_i(\boldsymbol{x}^k,{\boldsymbol{y}}^k)+ \sum\limits_{j=1}^r v^{k+1}_j \nabla h_j(\boldsymbol{x}^k,{\boldsymbol{y}}^k)=0,\label{eq:d33}
\end{eqnarray}
where $\bar{u}^{k}_i=u^k_i+2\rho_kg_i^+(\boldsymbol{x}^k,{\boldsymbol{y}}^k)\geq 0,i=1,2,\cdots,s$ and $v^{k+1}_j=v^k_j+2\rho_kh_j(\boldsymbol{x}^k,{\boldsymbol{y}}^k),j=1,2,\cdots,r$.
Let
 \begin{eqnarray*}
\gamma_k=1+\sum\limits_{i=1}^s \bar{u}^{k}_i+\sum\limits_{j=1}^r (\max\{v^{k+1}_j,0\}+\max\{-v^{k+1}_j,0\})>0.
\end{eqnarray*}
 Let $\eta^k=\frac{1}{\gamma_k}>0$, $\alpha_i^k=\frac{\bar{u}^{k}_i}{\gamma_k}\geq 0(i=1,2,\cdots,s)$;
$\mu_j^k=\frac{\max\{v^{k+1}_j,0\}}{\gamma_k}\geq 0(j=1,2,\cdots,r)$ and $\nu_j^k=\frac{\max\{-v^{k+1}_j,0\}}{\gamma_k}\geq 0(j=1,2,\cdots,r)$.
Then,
 \begin{eqnarray}
 \eta^k+\sum\limits_{i=1}^s \alpha_i^k+\sum\limits_{j=1}^r(\mu_j^k+\nu_j^k)=1.\label{eq:d34}
 \end{eqnarray}
Clearly, as $k\to\infty$, we have $\eta^k\rightarrow \eta>0, \alpha_i^k\rightarrow \alpha_i(i=1,2,\cdots,s), \mu_j^k\rightarrow \mu_j(j=1,2,\cdots,r)$ and $\nu_j^k\rightarrow \nu_j(j=1,2,\cdots,r)$.
By \eqref{eq:d33} and \eqref{eq:d34}, we have
 \begin{eqnarray}
\eta\nabla g(\boldsymbol{x}^*,{\boldsymbol{y}}^*)+\sum\limits_{i=1}^s\alpha_i\nabla g_i(\boldsymbol{x}^*,{\boldsymbol{y}}^*)+ \sum\limits_{j=1}^r( \mu_i -\nu_i)\nabla h_j(\boldsymbol{x}^*,{\boldsymbol{y}}^*)=0.\label{eq:d35}
\end{eqnarray}
By  \eqref{eq:d35} and Lemma 3, let $\beta^k_i=\mu_i^k-\nu_i^k\to \beta_i$ as $k\rightarrow +\infty$,  \eqref{eq:d31} and \eqref{eq:d32} hold.\\

Now, $((\boldsymbol{x}_1,\boldsymbol{y}_1),(\boldsymbol{x}_2,\boldsymbol{y}_2),\cdots,(\boldsymbol{x}_p,\boldsymbol{y}_p))$ is called a decomposition of $(\boldsymbol{x},\boldsymbol{y}))$ on $S$ if it satisfies the following conditions:

(i) $\boldsymbol{x}=(\boldsymbol{x}_{1}, \boldsymbol{x}_{2},\cdots,\boldsymbol{x}_{p})^\top\in S$, where $\boldsymbol{x}_j=(x_{j1},x_{j2},\cdots,x_{jp_j})\in R^{p_j}, j=1,2,\cdots,p$ and $ \sum\limits_{j=1}^p p_j=n$;

(ii) $\boldsymbol{y}=(\boldsymbol{y}_{1}, \boldsymbol{y}_{2},\cdots,\boldsymbol{y}_{p})^\top\in R^m$, where $\boldsymbol{y}_j=(y_{j1},y_{j2},\cdots,y_{jq_j})\in R^{q_j}, j=1,2,\cdots,p$ and $ \sum\limits_{j=1}^p q_j=m$;

(iii) $((\boldsymbol{x}_1,\boldsymbol{y}_1),(\boldsymbol{x}_2,\boldsymbol{y}_2),\cdots,(\boldsymbol{x}_p,\boldsymbol{y}_p))$ is a rearrangement of $(\boldsymbol{x},\boldsymbol{y})$;

(iv) there are no identical variables between $(\boldsymbol{x}_k,\boldsymbol{y}_k)$ and $(\boldsymbol{x}_j,\boldsymbol{y}_j)$ for any $j,k=1,2,\cdots,p, k\not=j)$.

 Let $(\boldsymbol{x}_j,\boldsymbol{y}_j\mid(\boldsymbol{x},\boldsymbol{y})):=((\boldsymbol{x}_1,\boldsymbol{y}_1),(\boldsymbol{x}_2,\boldsymbol{y}_2),\cdots,(\boldsymbol{x}_p,\boldsymbol{y}_p))$, where $(\boldsymbol{x}_j,\boldsymbol{y}_j)$ is a variable, i.e. all $(\boldsymbol{x}_k,\boldsymbol{y}_k)(k=1,2,\cdots,p, k\not=j)$ are fixed except  $(\boldsymbol{x}_j,\boldsymbol{y}_j)$.

 Let $R^{s_j}_+=\{\boldsymbol{\beta}_j\in R^{s_j}\mid \boldsymbol{\beta}_j\geq 0\},j=1,2,\cdots,p$. Let $ \boldsymbol{\alpha}_j\in R^{r_j},\boldsymbol{\beta}_j\in R^{s_j}_+,j=1,2,\cdots,p$ be Lagrange parameters and $\sigma>0$ be a penalty parameter with $(\boldsymbol{x},\boldsymbol{y})=((\boldsymbol{x}_1,\boldsymbol{y}_1),(\boldsymbol{x}_2,\boldsymbol{y}_2),\cdots,(\boldsymbol{x}_p,\boldsymbol{y}_p))$, $\boldsymbol{u}=(\boldsymbol{\alpha}_1,\boldsymbol{\alpha}_2,\cdots,\boldsymbol{\alpha}_p)$ and $\boldsymbol{v}=(\boldsymbol{\beta}_1, \boldsymbol{\beta}_2,\cdots,\boldsymbol{\beta}_p)$.

In order to solve (CNP), we reduce the scale problem of (CNP) and use the decomposition method to solve (CNP) with ALPF Algorithm. The augmented Lagrange penalty functions of all subproblems (CNP)$_j(\boldsymbol{\alpha}_j,\boldsymbol{ \beta}_j,\sigma)(j=1,2,\cdots,p)$ are defined by
\begin{eqnarray}
A_j(\boldsymbol{x}_i,\boldsymbol{y}_j\mid(\boldsymbol{x},\boldsymbol{y});\boldsymbol{\alpha}_j,\boldsymbol{\beta}_j, \sigma)&=&g(\boldsymbol{x}_i,\boldsymbol{y}_j|(\boldsymbol{x},\boldsymbol{y}))+\boldsymbol{\alpha}_j^\top \boldsymbol{h}_{j}(\boldsymbol{x}_{j},\boldsymbol{y}_{j}\mid(\boldsymbol{x},\boldsymbol{y})\nonumber \\ &&
+\boldsymbol{\beta}_j^\top \boldsymbol{g}_{j}(\boldsymbol{x}_{j},\boldsymbol{y}_{j}\mid(\boldsymbol{x},\boldsymbol{y})
+\frac{1}{2}\sigma \|\boldsymbol{h}_{j}(\boldsymbol{x}_{j},\boldsymbol{y}_{j}\mid(\boldsymbol{x},\boldsymbol{y})\|^2\nonumber \\ && +\frac{1}{2}\sigma \max\{\boldsymbol{g}_{j}(\boldsymbol{x}_{j},\boldsymbol{y}_{j}\mid(\boldsymbol{x},\boldsymbol{y}),0\}^2, \label{eq:d45}
\end{eqnarray}
where  $(\boldsymbol{x}_j,\boldsymbol{y}_j)$ is variable, i.e. all $(\boldsymbol{x}_k,\boldsymbol{y}_k)(k=1,2,\cdots,p, k\not=j)$ are fixed except for $(\boldsymbol{x}_j,\boldsymbol{y}_j)$. By \eqref{eq:d45}, for $j=1,2,\cdots,p$,  unconstraint optimization subproblems are defined by
\begin{eqnarray*}
\mbox{(CNP)}_j(\boldsymbol{\alpha}_j,\boldsymbol{\beta}_j,\sigma)\qquad & \min\; & A_j(\boldsymbol{x}_j,\boldsymbol{y}_j|(\boldsymbol{x},\boldsymbol{y});\boldsymbol{\alpha}_j,\boldsymbol{\beta}_j,\sigma) \\
& \mbox{s.t.}\; & (\boldsymbol{x}_j,\boldsymbol{y}_j)\in R^{p_j}\times R^{q_j}.
\end{eqnarray*}

For $j=1,2,\cdots,p$ and $n=ep$, the subproblems (CNP)$_j$ are solved by repeatedly using the ALPF Algorithm such that larger scale problem (CNP) may be solved.

Finally, a nonconvex optimization problem with convertible nonconvex function $f$ is defined by
 \begin{eqnarray*}
\mbox{(UOP)} &\min& f(\boldsymbol{x})\\
&s.t.& \boldsymbol{x}\in S,\\
\end{eqnarray*}
where $S$ is nonempty set in $R^n$. If $f$ is a differentiable convertible nonconvex function with $f=[g:g_1,g_2,\cdots, g_s; h_1,h_2,\cdots,h_r]$, then $A(\boldsymbol{x},\boldsymbol{y};\boldsymbol{u},\boldsymbol{v},\rho)$ is defined by  \eqref{eq:d26}. So , the following problem is solved by ALPF Algorithm.
\begin{eqnarray*}
\mbox{(CNP)}(\boldsymbol{u},\boldsymbol{v},\rho) &\min& A(\boldsymbol{x},\boldsymbol{y};\boldsymbol{u},\boldsymbol{v},\rho), \\
                        &s.t.&\boldsymbol{x}\in S, \boldsymbol{y}\times R^m.
\end{eqnarray*}
We have obtained numerical results of three examples with Matlab.

{\bf Example 7} A nonconvex optimization problem is  (Problem 1 in \cite{Zhang})
 \begin{eqnarray*}
\mbox{(EX7)} &\min& f(\boldsymbol{x})=2x_1^2-1.05x_1^4+\frac{1}{6}x_1^6-x_1x_2+x_2^2\\
&s.t.& |x_1|\leq 3,|x_2|\leq 3.
\end{eqnarray*}
An optimal solution to (EX7) is (0,0) in \cite{Zhang}. Let $(\boldsymbol{x},\boldsymbol{y})=(x_1,x_2,y_1,y_2,y_3)$. A convertible nonconvex form of $f$ is defined by
 \begin{eqnarray*}
f(\boldsymbol{x})=[2x_1^2-1.05y_1+\frac{1}{6}y_2+0.5(x_1-x_2)^2-0.5y_3+x_2^2:\\
-x_1-3\leq 0,x_1-3\leq 0,-x_2-3\leq 0,-x_4-3\leq 0;\\
x_1^4-y_1=0,x_2^6-y_2=0,x_1^2+x_2^2-y_3=0].
\end{eqnarray*}
In ALPF Algorithm, the starting parameters $\epsilon=10^{-6},\rho_{1}=10,N=100,\boldsymbol{u}^0=(0, 0, 0)$ and
$(\boldsymbol{x}^0,\boldsymbol{y}^0)=(2, 2, 2,2,2)$ are taken.  At the 3th step, an approximate solution
$\boldsymbol{x}^3=(0.000000,0.000000)$ to (EX7) is obtained.

{\bf Example 8}  A nonconvex nonsmooth optimization problem is (Problem 5 in \cite{Bagirov})
 \begin{eqnarray*}
\mbox{(EX8)} &\min& f_n(\boldsymbol{x})=n\max\{|x_i|: i=1,2,\cdots,n\}-\sum\limits_{i=1}^n |x_i| \\
&s.t.& \boldsymbol{x}\in R^n.
\end{eqnarray*}
An optimal solution to (EX8) is $\boldsymbol{x}^*=(\pm\alpha,\pm\alpha,\cdots,\pm\alpha)^\top$ in \cite{Bagirov} with $f(\boldsymbol{x}^*)=0$ for $\alpha \in R^1$ in \cite{Bagirov}. Let $\boldsymbol{x}\in R^n,\boldsymbol{y}\in R^{2n+1}$. A convertible nonconvex form of $f$ is defined by
 \begin{eqnarray*}
g(\boldsymbol{x},\boldsymbol{y})&=&n y_{2n+1}-\sum\limits_{i=1}^n y_i \\
h_i(\boldsymbol{x},\boldsymbol{y})&=& y_i^2-y_{i+n}=0,\ i=1,2,\cdots,n,\\
h_{n+i}(\boldsymbol{x},\boldsymbol{y})&=&x_i^2-y_{i+n}=0,\ i=1,2,\cdots,n,\\
g_{i}(\boldsymbol{x},\boldsymbol{y})&=&-y_i\leq 0,\ i=1,2,\cdots,n,\\
g_{n+i}(\boldsymbol{x},\boldsymbol{y})&=&y_i-y_{2n+1}\leq 0, \ i=1,2,\cdots,n.
\end{eqnarray*}
In  ALPF Algorithm, the starting parameters $\epsilon=10^{-6},\rho_{1}=10,N=100,(\boldsymbol{u}^0,\boldsymbol{v}^0)=(0, 0, \cdots,0 )$ and
$(\boldsymbol{x}^0,\boldsymbol{y}^0)=(1, 2, 3,\cdots,3n+1)$ are taken.  For $n=5$, at the 4th step, an approximate solution
$\boldsymbol{x}^4=(3.4366,-3.4366,-3.4366, 3.4366,3.4366)$ to (EX8) is obtained.

We can also use the penalty function method to solve (CNP). Define a penalty function of (CNP) as follows
\begin{eqnarray*}
\mbox{(CNP)}(\rho)&\min& F(\boldsymbol{x},\boldsymbol{y};\rho)=g(\boldsymbol{x},\boldsymbol{y})+\rho \sum\limits_{i=1}^s g_i^+(\boldsymbol{x},\boldsymbol{y})^2+\rho \sum\limits_{j=1}^r h_j(\boldsymbol{x},\boldsymbol{y})^2,\nonumber\\
&s.t.& \ \ (\boldsymbol{x},\boldsymbol{y})\in R^n\times R^m,
\end{eqnarray*}
where $\rho>0$ is a penalty parameter and $g_i^+(\boldsymbol{x},\boldsymbol{y})=\max\{g_i(\boldsymbol{x},\boldsymbol{y}),0\}$.
By using Matlab, the starting parameters $\epsilon=10^{-6},\rho_{1}=10,N=100$ and
$(\boldsymbol{x}^0,\boldsymbol{y}^0)=(1, 2, 3,\cdots,3n+1)$ are taken in (CNP)$(\rho)$ for (EX8).  For $n=10$, at the 4th step, an approximate solution $\boldsymbol{x}^4=(-3.4015,3.4015,$ $3.4015,3.4015,3.4015,3.4015,3.4015,$ $3.4015,3.4015,-3.4015)$ to (EX8) is obtained.

For $n=500$ and $p=100$, the subproblems (CNP)$_j$ are solved by repeatedly $100$ times using decomposed ALPF Algorithm. An approximate solution $\boldsymbol{x}^4=(3.3484,3.3484,3.3484,\cdots,-3.3484)$ is obtained.

{\bf Example 9} A especial case in Example 4 (in \cite{Chen1}) is defined by
 \begin{eqnarray*}
\mbox{(EX9)} &\min& f_n(\boldsymbol{x})=(\sum\limits_{i=1}^n ix_i-2n)^2+\lambda\|\boldsymbol{x}\|_0 \\
&s.t.& \boldsymbol{x}\in R^n.
\end{eqnarray*}
 For $(x_i,y_i,y_{i+n})$, $i=1,2,\cdots,n$, a convertible nonconvex form of $f_n(\boldsymbol{x})$  is defined by
\begin{eqnarray*}
[(\sum\limits_{i=1}^n ix_i-2n)^2+\lambda\sum\limits_{i=1}^n y_{i}^2 :
((x_i+y_i-1)^2-y_{i+n},\\
x_i^2+(y_i-1)^2-y_{i+n},y_{i}^2-y_i)=0, i=1,2,\cdots,n].
\end{eqnarray*}
Let $e^k(\boldsymbol{x}^k,\boldsymbol{y}^k)=\|\boldsymbol{g}^+(\boldsymbol{x}^k,\boldsymbol{y}^k)\|+\|\boldsymbol{h}(\boldsymbol{x}^k,\boldsymbol{y}^k)\|$. Let $n=10$, $\epsilon=0.000001,\rho_{1}=10,N=10, (\boldsymbol{x}^0,\boldsymbol{y}^0)\in R^n\times R^m$, $\boldsymbol{v}^1\in R^r$,$(\boldsymbol{x}^0,\boldsymbol{y}^0)=(0, 0, 0,\cdots,0)$, $\boldsymbol{u}^1=\boldsymbol{g}^+(\boldsymbol{x}^0,\boldsymbol{y}^0)$,$\boldsymbol{v}^{1}=0$. The numerical results of (EX9) are obtained by ALPF Algorithm at $\lambda_1=1$ and $\lambda_1=10$ in Table 1 and Table 2. Table 1 and Table 2 show that approximate sparse solution may be obtained by ALPF Algorithm when the penalty parameter increases.

\begin{center}
{\small Table 1.  The numerical results of (EX9) are obtained by ALPF Algorithm at $\lambda_1=1$.   \\
\begin{tabular}{cccccc}
 \hline
  $k$  &$\rho_k$ &$\boldsymbol{x}^k$  &$f(\boldsymbol{x}^k)$&$\|\boldsymbol{x}^k\|_0$&$e^k(\boldsymbol{x}^k,\boldsymbol{y}^k)$\\
   \hline
1 & 5 &(0.0021,0.0044,0.0059,0.0072,0.0088,& 1.9135 & 1.9131  &0.1372\\
&&0.0119,0.0129,0.0202,-0.4906,2.3972)&&&\\
2 & 50 &(-0.0003,-0.0005,-0.0005,-0.0007,-0.0008,& 2.0066 & 2.0066  &0.0063\\
&&      -0.0011,-0.0012,-0.0018,-1.1625,3.0500)\\
3 & 500 &(0.0000,0.0000,-0.0000,-0.0000,-0.0000,& 2.0002 & 2.0002  &0.0001\\
&&-0.0000,-0.0000,-0.0000,-0.9386,2.8448)&&&\\
 \hline
\end{tabular}
}
\end{center}

\begin{center}
{\small Table 2.  The numerical results of (EX9) are obtained by ALPF Algorithm at $\lambda_1=10$.   \\
\begin{tabular}{cccccc}
 \hline
  $k$  &$\rho_k$ &$\boldsymbol{x}^k$  &$f(\boldsymbol{x}^k)$&$\|\boldsymbol{x}^k\|_0$&$e^k(\boldsymbol{x}^k,\boldsymbol{y}^k)$\\
   \hline
1 & 5 &(0.0506,0.1013,0.1519,0.2025,0.2532, & 0.2564 & 0.0000  &1.4050\\
 &  & 0.3038,0.3544,0.4051,0.4557,0.5063) &   &    & \\
2 & 50 &(-0.0051,-0.0101,-0.0152,-0.0203,-0.0253,& 10.0000 & 1.0000  &0.1209\\
  &    & -0.0304,-0.0354,-0.0405,-0.0456,2.1443) &&&\\
3 & 500 &(0.0000,0.0000,0.0000,0.0000,0.0000,& 10.0000 & 1.0000  &0.0000\\
&         &0.0000,0.0000,0.0000,0.0000,2.0000) &&&\\
\hline
\end{tabular}
}
\end{center}

By  ALPF Algorithm,let the starting parameters $\epsilon=10^{-4},\sigma_{1}=5,N=10, \boldsymbol{\alpha}_j^0=(0, 0, \cdots,0)$ and $\boldsymbol{\beta}_j^0=(0,0,0,\cdots,0)$ be taken. The subproblems (CNP)$_j(\boldsymbol{\alpha}_j, \boldsymbol{\beta}_j,\rho)(j=1,2,\cdots,p)$ are solved by repeatedly $p$ times using the ALPF Algorithm. Let $n=ep$. When $e=5,p=1,2,6,10,20,$ $100,200$,$\lambda=1,10$, values of
0-norm $\|\boldsymbol{x}^k\|_0$ are obtained by ALPF Algorithm in Table 3. Numerical results show that values of
0-norm $\|\boldsymbol{x}^k\|_0$ decrease in Table 3 when $\lambda$ increases.

\begin{center}
 Table 3.  The values of 0-norm $\|\boldsymbol{x}^k\|_0$ of (EX9) are obtained by ALPF Algorithm.   \\
\begin{tabular}{ccc}
 \hline
$n=5p$& $\lambda=1$ & $\lambda=10$ \\
\hline
5  & 1 & 1 \\
10 & 2 & 2 \\
30 &15  & 10 \\
50 & 22 &17  \\
100 &54 &37  \\
500 &325 &192 \\
1000 &611&382 \\
\hline
\end{tabular}
\end{center}

Hence, the above examples  illustrate that it is efficient to solve an approximate optimal solution to (UOP) by using the ALPF algorithm with Maltlab, to avoids using subdifferentiation or smoothing techniques.

\section{Conclusion}

In this paper, a new concept-(exact) convertible nonconvex function is proposed. It involves many nonconvex and nonsmooth functions, even discontinuous nonconvex functions. The convertible nonconvex function is a set of convex functions which map the known nonconvex function to the constrained equalities and inequalities in high dimensional space. These constraint functions are convex. The sufficient condition of global optimal solution is proved, which is equivalent to KKT condition. The strong duality theorem of Lagrange function of convertible nonconvex function is proved. Therefore,  an augmented Lagrangian penalty function algorithm is proposed and its convergence is proved. If there is a global optimal solution to the optimization problem of convertible nonconvex function, then the algorithm may be find a global optimal solution theoretically. Numerical results show that it is efficient to solve (UOP) by  ALPF Algorithm avoids using subdifferentiation or smoothing techniques, allowing a direct use of some gradient search algorithms, such as gradient descent algorithm, Newton algorithm and so on.

This paper provides a new idea for solving nonconvex or nonsmooth optimization problems in theory. It has important potential value for solving nonconvex or nonsmooth optimization problems in many application fields.


\bibliographystyle{siamplain}

\end{document}